\documentclass[11pt,reqno,letterpaper]{amsart}

\usepackage[T1]{fontenc}
\usepackage{lmodern}
\usepackage{amsmath,amssymb,amsthm,mathtools}
\usepackage{enumitem,mathrsfs}
\usepackage{microtype}
\usepackage[hidelinks,bookmarksnumbered]{hyperref}
\allowdisplaybreaks[2]
\numberwithin{equation}{section}
\makeatletter
\def\@settitle{\begin{center}%
  \Large\bfseries
  \baselineskip18\p@\relax
  \@title
  \end{center}}
\makeatother

\newtheorem{theorem}{Theorem}[section]
\newtheorem{proposition}[theorem]{Proposition}
\newtheorem{lemma}[theorem]{Lemma}
\newtheorem{corollary}[theorem]{Corollary}
\theoremstyle{definition}

\newtheorem{definition}[theorem]{Definition}
\newtheorem{question}[theorem]{Question}

\DeclareMathOperator{\Ric}{Ric}
\DeclareMathOperator{\biRic}{biRic}
\DeclareMathOperator{\Reg}{Reg}
\DeclareMathOperator{\Sing}{Sing}
\DeclareMathOperator{\Vol}{Vol}
\DeclareMathOperator{\spt}{spt}
\DeclareMathOperator*{\essinf}{ess\,inf}
\DeclareSymbolFont{averageintegrals}{U}{esint}{m}{n}
\DeclareMathSymbol{\fintop}{\mathop}{averageintegrals}{'037}
\newcommand{\fint}{\fintop\nolimits}

\title[Bi-Ricci curvature and volume in all dimensions]
{Sharp volume estimates under\\ positive bi-Ricci curvature in all dimensions}
\author{Jintian Zhu}
\address[Jintian Zhu]{Institute for Theoretical Sciences, Westlake
University, 600 Dunyu Road, 310030, Hangzhou, Zhejiang,
People's Republic of China}
\email{zhujintian@westlake.edu.cn}

\subjclass[2020]{Primary 53C21; Secondary 53C20, 49Q05}
\keywords{Bi-Ricci curvature, volume growth, prescribed mean curvature,
spectral Bishop inequality, minimal hypersurface singularities}

\begin{document}

\begin{abstract}
Let $(M^n,g)$, $n\geq3$, be a connected complete noncompact Riemannian
manifold with nonnegative Ricci curvature. Suppose that its bi-Ricci
curvature is at least $n-2$ outside a compact set. We prove
$$
 \Vol B_R(p)\leq2\omega_{n-1}R
$$
for every $p\in M$ and $R>0$, where $\omega_{n-1}$ denotes the area
of the unit round $(n-1)$-sphere. We also prove
$$\lim_{R\to\infty}\frac{\Vol B_R(p)}{R}\leq e(M)\omega_{n-1},$$ where $e(M)$ denotes the number
of ends of $M$. Both estimates above are sharp. The proof combines almost
minimal prescribed-mean-curvature boundaries in noncompact manifolds of nonnegative Ricci curvature
 with an extension of the
spectral Bishop-Gromov volume comparison by Antonelli and Xu \cite{AntonelliXu} to AM-PI spaces where
mild singularities are allowed.
\end{abstract}

\maketitle

\section{Introduction}
\label{sec:introduction}

In the study of scalar curvature, Gromov \cite[Section~2.A(b)]{Gromov} asked whether a complete noncompact
Riemannian manifold $(M^n,g)$, $n\geq3$, with $\Ric_M\geq0$ and
$\operatorname{Scal}_M\geq2$ satisfies
$$
 \Vol B_R(p)\leq C(n)R^{n-2}
 \,\,\text{for every }p\in M\text{ and }R>0.
$$
Munteanu and Wang \cite{MunteanuWang} proved the conjecture in dimension
three. Chodosh, Li, and Stryker \cite{ChodoshLiStryker} subsequently
gave a proof of linear volume growth using $\mu$-bubbles, and Wang
\cite{WangVolume} refined this approach to obtain a universal constant.
Wei, Xu, and Zhang \cite{WeiXuZhang} obtained the sharp asymptotic
volume growth bound in dimension three and the corresponding rigidity.
Three recent preprints independently by Ge \cite{Ge},
Kong and Zhu \cite{KongZhu}, and Antonelli \cite{AntonelliGrowth}
propose proofs of the conjecture in all dimensions, while the sharp 
constant $C(n)$ remains open.

Motivated by the work of Antonelli and Xu \cite{AntonelliXu} and
Antonelli \cite{AntonelliGrowth}, we study sharp volume bounds under
positive bi-Ricci curvature. Introduced by Shen and Ye \cite{ShenYe},
the bi-Ricci curvature is defined for an orthonormal pair $v,w\in T_xM$ by
$$
 \biRic_M(v,w)
 =\Ric_M(v,v)+\Ric_M(w,w)-\operatorname{Rm}_M(v,w,w,v).
$$
It can be viewed as a generalization of the scalar curvature, since in dimension three we have
$
 \biRic_M(v,w)=\operatorname{Scal}_M/2
$.
For convenience, we write 
$$\biRic_M\geq k$$ when we have 
$$\biRic_M(v,w)\geq k$$ for every
orthonormal pair $(v,w)$ in the region under consideration.

Suppose that
 we have $\Ric_M\geq0$ on $M$ and also $\biRic_M\geq n-2$ outside a
compact set. In dimensions up to five, Antonelli and Xu
\cite{AntonelliXu} proved linear volume growth. When $M$ has one end,
they obtained the sharp asymptotic bound under the additional assumption
$\inf_{x\in M}\Vol B_1(x)>0$. They also posed the following question
\cite[Question~1]{AntonelliXu}.

\begin{question}\label{ques:antonelli-xu}
Let $(M^n,g)$, $n\geq6$, be a connected complete noncompact Riemannian
manifold with $\Ric_M\geq0$ on $M$ and $\biRic_M\geq n-2$ outside a
compact set. Does $M$ have linear volume growth? If $M$ has one end,
does its isoperimetric profile $I_M$ satisfy
$$
 \lim_{v\to\infty}I_M(v)\leq\omega_{n-1}?
$$
\end{question}

Zhou and the author \cite{ZhouZhu} later removed the
noncollapsing assumption and proved the sharp global and asymptotic
volume bounds in dimensions three through seven, answering this question
in dimensions six and seven. Recently, Antonelli's theorem on intermediate
curvatures \cite{AntonelliGrowth} yields linear volume growth with a
dimensional constant in every dimension when the bi-Ricci lower bound
holds on all of $M$. 

In this paper, we establish the sharp global and asymptotic volume
bounds in all dimensions assuming the bi-Ricci lower bound only outside
a compact set, extending the results of Zhou and the author \cite{ZhouZhu}.

Write $e(M)$ for the number of ends of $M$.
Our main theorem can be stated as follows.

\begin{theorem}\label{thm:main}
Let $(M^n,g)$, $n\geq3$, be a connected complete noncompact Riemannian manifold. Suppose that there is a compact set $K\subset M$ such
that
$$
 \begin{gathered}
 \Ric_M\geq0\,\,\text{on }M\\
 \text{and}\,\,\biRic_M\geq n-2\,\,\text{on }M\setminus K.
 \end{gathered}
$$
Then, for every $p\in M$ and $R>0$, we have
\begin{equation}\label{eq:globalsharp}
 \Vol B_R(p)\leq2\omega_{n-1}R.
\end{equation}
Moreover, the limit
\begin{equation}\label{eq:slope}
 \Theta^1(\infty)=\lim_{R\to\infty}\frac{\Vol B_R(p)}{R}
\end{equation}
exists, is independent of $p$, and satisfies
$$
 \Theta^1(\infty)\leq e(M)\omega_{n-1}.
$$
All these estimates are sharp.
\end{theorem}

The theorem also yields the one-ended isoperimetric bound at infinity.

\begin{corollary}\label{cor:one-ended-isoperimetry}
Let $(M^n,g)$, $n\geq3$, be a connected complete noncompact Riemannian
manifold with one end. Suppose that $\Ric_M\geq0$ on $M$ and
$\biRic_M\geq n-2$ outside a compact set. Then its isoperimetric
profile satisfies
$$
 \lim_{v\to\infty}I_M(v)\leq\omega_{n-1}.
$$
\end{corollary}

This answers Question~\ref{ques:antonelli-xu} in all dimensions.

The proof of Theorem~\ref{thm:main} applies the \emph{infinity-based} volume comparison approach 
of Zhou and the author \cite{ZhouZhu}. The idea is to enclose arbitrarily large balls by
separating hypersurfaces of area at most $\omega_{n-1}+o(1)$.
The sharp volume growth bound then follows from radial comparison along
a geodesic ray, with the associated Busemann function controlling the
annular width.

As before, we construct the desired separating hypersurfaces as components of locally minimizing
prescribed-mean-curvature boundaries.
In dimensions $3\leq n\leq7$, Zhou and the author \cite{ZhouZhu}
choose the prescribed mean curvature $h$ on a wide annular region so that
\begin{equation}\label{eq:annular-prescription}
 c_nh^2-|\nabla h|\geq-\varepsilon,
\end{equation}
for a dimensional constant $c_n>0$, where $\varepsilon$ is a small positive
constant satisfying $\varepsilon\to0$ as the width tends $+\infty$.
For the separating hypersurface $\Sigma$ with $\mathbf H=h$,
the Gauss equation and the assumption $\Ric_M\geq0$ give
$$
 \Ric_M(\nu,\nu)+|A|^2
 \geq\beta_n\bigl(\biRic_M(\nu,e)-\Ric_\Sigma(e,e)\bigr)+c_nh^2
$$
for every unit tangent vector $e$, where $\nu$ and $A$ denote the unit
normal and second fundamental form of $\Sigma$. Together with the
stability inequality for the prescribed mean curvature functional,
these inequalities give a spectral lower bound for $\Ric_\Sigma$
and hence the desired volume estimate. 

In dimensions greater than seven,
however, the spectral Bishop--Gromov volume comparison theorem requires
$\beta_n\geq(n-3)/(n-2)$, which forces $c_n<0$, and so
\eqref{eq:annular-prescription} cannot hold if the prescribed function $h$
has the blow-up behavior:
$h\to+\infty$ at the inner boundary and $h\to-\infty$ at the outer
boundary. Additional difficulty
 also comes from the fact that the minimizing
boundaries may also be singular.

The strategy here is to construct the separating hypersurface by prescribing
mean curvature function $h$ to have arbitrarily small $L^\infty$-norm. 
The key observation is that we are able to construct a Caccioppoli exhaustion, whose boundary mean curvature
tends uniformly to zero, in a complete noncompact Riemannian manifold
of nonnegative Ricci curvature. We then minimize a prescribed mean curvature
functional between two suitable sets of this exhaustion. Their perimeter
comparison inequalities allow us to choose $h_j$ linear in a smooth
approximation of the radial distance. Along the minimizing boundary, we have
$$
 |h_j|=O(j^{-1})\,\,\text{and}\,\,|\nabla h_j|=O(j^{-2}).
$$
These bounds ensure that \eqref{eq:annular-prescription} holds with
$\varepsilon=O(j^{-2})$ even when $c_n<0$. The stability inequality
for the minimizing boundary and the Gauss equation give a spectral Ricci lower bound with an $O(j^{-2})$ error
on the regular parts of the resulting separating hypersurfaces.

To obtain the sharp area bound for these possibly singular hypersurfaces,
we extend the weighted isoperimetric
profile argument of Antonelli and Xu \cite{AntonelliXu}. In order to establish
the ordinary differential
inequality for the weighted isoperimetric profile, 
the second
variation requires a cutoff function near the singular set, where the weight function
may be unbounded. A perimeter comparison
and an energy estimate for the logarithm of the weight function give an
upper bound for the weighted area in small balls. Combined with the small Hausdorff dimension
of the singular set, we can construct cutoff functions, whose
weighted Dirichlet energy on the minimizing boundary tends to zero.
The isoperimetric profile argument therefore works smoothly, yielding the sharp
area bound as desired.

\subsection*{Arrangements of the paper}

The rest of this paper is arranged as follows.
In Section~\ref{sec:exhaustion}, we construct a Caccioppoli exhaustion
whose boundary mean curvature tends uniformly to zero.
In Section~\ref{sec:spectral-comparison}, we establish a sharp intrinsic
spectral volume comparison theorem for compact AM--PI spaces.
Section~\ref{sec:volume} combines these results with radial volume
comparison to prove Theorem~\ref{thm:main}.
Appendix~\ref{sec:ampi} collects the necessary properties of compact
AM--PI spaces and shows that compact locally minimizing
prescribed-mean-curvature boundaries belong to this class.

\section*{Acknowledgements}
J. Zhu is partially supported by National Key R\&D Program of China
2023YFA1009900, NSFC Grant No.\,12401072, Zhejiang Provincial Natural
Science Foundation of China under Grant No.\,LQKWL26A0101 and the
start-up fund from Westlake University.

\section*{Declaration of AI usage}
The mathematical ideas are due to the author. The author used
ChatGPT (Astra 6.2) to help complete the area comparison argument establishing
boundedness of the minimizer in Lemma~\ref{lem:direct-boundary-properties}
and the analysis in Lemma~\ref{lem:ampi-weight}. It was also used to
check references and revise the exposition for clarity. The author takes
full responsibility for the mathematical content.

\section{Caccioppoli exhaustion with almost minimal boundary}
\label{sec:exhaustion}

In this section, we construct a Caccioppoli exhaustion with almost
minimal boundary in any 
connected, complete, noncompact Riemannian manifold $(M^n,g)$ of dimension
$n\geq2$ with nonnegative Ricci curvature.

We will use the following convention of notations. 
For a smooth boundary $\Gamma=\partial\Omega$, we use $\nu$ to denote the outward unit
normal and write $$\mathbf H=\operatorname{div}_\Gamma\nu.$$
For a Caccioppoli set $E\subset M$,
that is, a measurable set with $\chi_E\in BV_{\rm loc}(M)$, we write
$$
 \mathcal P(E;U)=|D\chi_E|(U)
 \,\,\text{and}\,\,\mathcal P(E)=\mathcal P(E;M)
$$
for its perimeter in an open set $U\subset M$ and its total perimeter,
as in \cite{WangSmoothing}.

Our main purpose is the following proposition.

\begin{proposition}\label{prop:caccioppoli-exhaustion}
There exist precompact open Caccioppoli sets $\{E_i\}_{i=1}^\infty$
such that
$$
 \overline{E_i}\subset E_{i+1}
 \,\,\text{and}\,\,\bigcup_{i=1}^\infty E_i=M.
$$
Each boundary 
$$\Gamma_i=\partial E_i=\spt|D\chi_{E_i}|$$ is
smooth away from a closed singular set of Hausdorff dimension at most $n-8$,
which is empty when $n\leq7$. The associated multiplicity-one varifolds
have generalized mean curvature vectors $\mathbf H_i$ satisfying
$$
 \lim_{i\to\infty}
 \|\mathbf H_i\|_{L^\infty(\Gamma_i)}=0.
$$
\end{proposition}

To construct the exhaustion, it suffices to enclose any given large ball by a
precompact Caccioppoli set with a quantitative bound on its boundary
mean curvature.

\begin{lemma}\label{lem:prescribed-boundaries}
For every $o\in M$ and $r>0$, there exists an open Caccioppoli set $E$
such that
$$
 \overline{B_r(o)}\subset E\subset B_{2r}(o).
$$
For every open set $U\Subset M\setminus\overline{B_r(o)}$ and every
Caccioppoli set $F$ with $E\Delta F\Subset U$, we have
\begin{equation}\label{eq:exhaustion-almost-minimality}
 \mathcal P(E;U)\leq\mathcal P(F;U)
 +\frac{3n}{r}\Vol(E\Delta F).
\end{equation}
Its boundary $\Gamma=\partial E=\spt|D\chi_E|$ is smooth away from a
closed singular set of Hausdorff dimension at most $n-8$, which is
empty when $n\leq7$. The associated multiplicity-one varifold has
generalized mean curvature vector $\mathbf H$ satisfying
\begin{equation}\label{eq:direct-mean-curvature}
 \|\mathbf H\|_{L^\infty(\Gamma)}\leq\frac{3n}{r}.
\end{equation}
\end{lemma}

We apply the \emph{$\mu$-bubble method with an inner obstacle} as in \cite[Section~2.1]{BiZhu}.

Fix $o\in M$ and $r>0$. Choose a smooth cutoff function 
$\eta:M\to[0,1]$
supported in $B_{3r/2}(o)$ and equal to one on a neighborhood of
$\overline{B_r(o)}$. Set
$$
 h=\frac{3n}{r}(2\eta-1).
$$
Thus $h=3n/r$ near $\overline{B_r(o)}$, $h=-3n/r$ outside
$B_{3r/2}(o)$, and $|h|\leq3n/r$ on $M$.
Let $\mathcal C$ be the class of sets with finite volume and finite
perimeter that contain $B_r(o)$ modulo null sets. On $\mathcal C$,
consider the functional
$$
 \mathcal J(F)=\mathcal P(F)-\int_Fh\,dV.
$$

\begin{lemma}\label{lem:direct-existence}
The functional $\mathcal J$ attains its minimum on $\mathcal C$.
\end{lemma}

\begin{proof}
For $F\in\mathcal C$, we have
$$
\begin{aligned}
 \mathcal J(F)
 &=\mathcal P(F)-\int_Fh\,dV\\
 &=\mathcal P(F)+\frac{3n}{r}\Vol(F)-\frac{6n}{r}\int_F\eta\,dV.
\end{aligned}
$$
By the coarea formula for the distance function $d(o,\cdot)$, the ball
$B_s(o)$ has finite perimeter for almost every $s\in(r,2r)$.
Each such ball belongs to $\mathcal C$, yielding that $\mathcal C$ is non-empty. Since $0\leq\eta\leq1$ and $\eta$ has compact
support, $\mathcal J$ is bounded below. Choose a minimizing sequence
$\{F_i\}_{i=1}^\infty\subset\mathcal C$, then we have a uniform bound on
$\mathcal P(F_i)+(3n/r)\Vol(F_i)$.

By the local compactness theorem for sets of finite perimeter, a
subsequence converges in $L^1_{\rm loc}(M)$ to a set $E$.
Lower semicontinuity of perimeter and volume gives
$$\mathcal P(E)+\Vol(E)<\infty.$$ Since each $F_i$ contains $B_r(o)$
modulo null sets, the same holds for $E$. Thus we have $E\in\mathcal C$.
Since $\eta$ is bounded and compactly supported, local $L^1$ convergence
gives $$\int_{F_i}\eta\,dV\to\int_E\eta\,dV$$ along this subsequence.
It follows that $$\mathcal J(E)\leq\liminf_{i\to\infty}\mathcal J(F_i)
=\inf_{\mathcal C}\mathcal J,$$ so $E$ is a minimizer.
\end{proof}

To locate the minimizer, we first estimate the change in perimeter
when a geodesic ball is adjoined to or removed from a set.

\begin{lemma}\label{lem:ball-perimeter-comparison}
For every set $F$ of finite perimeter, every $p\in M$, and almost every
$s>0$, we have
$$
\begin{aligned}
 \mathcal P(F\setminus B_s(p))-\mathcal P(F)
 &\leq\frac ns\Vol(F\cap B_s(p)),\\
 \mathcal P(F\cup B_s(p))-\mathcal P(F)
 &\leq\frac ns\Vol(B_s(p)\setminus F).
\end{aligned}
$$
\end{lemma}

\begin{proof}
Write $d_p=d(p,\cdot)$ and $v(s)=\Vol(F\cap B_s(p))$.
The coarea formula for $d_p$ shows that $v$ is locally absolutely
continuous. By the slicing formula for sets of finite perimeter,
for almost every $s>0$ we have
$$
 \mathcal P(F\setminus B_s(p))-\mathcal P(F)
 =v'(s)-\mathcal P(F;B_s(p)).
$$
It remains to bound $v'(s)$.

Fix such an $s$ and let $0<\varepsilon<s$. Let $\zeta$ equal one on
$B_{s-\varepsilon}(p)$ and zero outside $B_s(p)$, with linear
interpolation in $d_p$ on the intervening annulus.
By the Laplacian comparison theorem, $\Delta(d_p^2/2)\leq n$ in the
sense of distributions. Testing with $\zeta u$ for a smooth
nonnegative function $u$ and using $|\nabla d_p|=1$ almost
everywhere, we obtain
$$
 \frac{s-\varepsilon}{\varepsilon}
 \int_{B_s(p)\setminus B_{s-\varepsilon}(p)}u\,dV
 \leq n\int_{B_s(p)}u\,dV+s\int_M\zeta|\nabla u|\,dV.
$$
Approximate $\chi_F$ in $L^1_{\mathrm{loc}}(M)$ by smooth
nonnegative functions whose variation measures converge weakly to
$|D\chi_F|$. Since $\zeta$ is continuous and
$0\leq\zeta\leq\chi_{B_s(p)}$, passing to the limit gives
$$
 (s-\varepsilon)\frac{v(s)-v(s-\varepsilon)}{\varepsilon}
 \leq n v(s)+s\mathcal P(F;B_s(p)).
$$
Letting $\varepsilon\downarrow0$ and substituting the bound for $v'(s)$
into the slicing formula proves the first inequality.
Applying it to $M\setminus F$ gives the second inequality, since complementary
sets have the same perimeter.
\end{proof}

We now apply Lemma~\ref{lem:ball-perimeter-comparison} to locate
every minimizer of $\mathcal J$ in $\mathcal C$.

\begin{lemma}\label{lem:direct-boundary-properties}
Every minimizer of $\mathcal J$ in $\mathcal C$ has an open
representative $E$ such that
$$
 \overline{B_r(o)}\subset E\subset B_{2r}(o).
$$
For every open set $U\Subset M\setminus\overline{B_r(o)}$ and every
Caccioppoli set $F$ with $E\Delta F\Subset U$, the perimeter comparison
\eqref{eq:exhaustion-almost-minimality} holds.
\end{lemma}

\begin{proof}
Let $E$ be a minimizer and let $\Omega$ be its essential interior, namely
$$
 \Omega=\bigl\{x\in M:\Vol(B_\rho(x)\setminus E)=0
 \,\,\text{for some }\rho>0\bigr\}.
$$
Then $\Omega$ is open and contains $B_r(o)$.
Let $U\Subset M\setminus\overline{B_r(o)}$ be open, and let $F$ be a
Caccioppoli set agreeing with $E$ almost everywhere outside a compact
subset of $U$. Clearly, we have $F\in\mathcal C$, and by the minimality of $E$,
\begin{equation}\label{eq:minimizer-perimeter-difference}
\begin{aligned}
 \mathcal P(E;U)-\mathcal P(F;U)
 &\leq\int_U(\chi_E-\chi_F)h\,dV\\
 &\leq\|h\|_{L^\infty(M)}\Vol(E\Delta F).
\end{aligned}
\end{equation}
Since $|h|\leq3n/r$, this proves \eqref{eq:exhaustion-almost-minimality}.

Suppose $U\Subset B_s(x)\Subset M\setminus\overline{B_r(o)}$.
Since $\chi_E-\chi_F$ has compact support in $U$, there exists a sequence
$\{u_k\}_{k=1}^\infty\subset C_c^\infty(U)$ converging strictly to
$\chi_E-\chi_F$ in $BV(U)$.
Since $\Ric_M\geq0$, we may apply the $L^1$ Poincar\'e inequality
\cite[Theorem~1.1]{LiSchoen} on $B_s(x)$ to each $u_k$.
Letting $k\to\infty$ gives
\begin{equation}\label{eq:symmetric-difference-volume}
 \Vol(E\Delta F)\leq Cs\bigl(\mathcal P(E;U)+\mathcal P(F;U)\bigr),
\end{equation}
where $C$ depends only on $n$. If $Cs\|h\|_{L^\infty(M)}\leq1/2$,
\eqref{eq:minimizer-perimeter-difference} and \eqref{eq:symmetric-difference-volume}
give $$\mathcal P(E;U)\leq3\mathcal P(F;U).$$
This yields that $E$ is locally quasiminimizing for perimeter outside
$\overline{B_r(o)}$ in the sense of
\cite[Definition~3.1]{KinnunenKorteLorentShanmugalingam}.
Hence $E$ agrees almost everywhere with its essential interior $\Omega$ in
$M\setminus\overline{B_r(o)}$; see
\cite[Theorem~4.2 and Remark~5.4]{KinnunenKorteLorentShanmugalingam} for instance.
Together with $\Vol(\partial B_r(o))=0$, this gives
$$\Vol(E\Delta\Omega)=0.$$ Thus $\Omega$ is an open representative of $E$.
For the rest of the proof, we identify $E$ with $\Omega$.

Choose $s>r$ close enough to $r$ that $h=3n/r$ on $B_s(o)$.
For almost every such $s$, Lemma~\ref{lem:ball-perimeter-comparison}
and the minimality of $E$ give
$$
 0\leq\mathcal J(E\cup B_s(o))-\mathcal J(E)
 \leq\left(\frac ns-\frac{3n}{r}\right)\Vol(B_s(o)\setminus E).
$$
Since $n/s<3n/r$, the ball $B_s(o)$ is contained in $E$ modulo null
sets. As $E$ is its essential interior, we obtain $\overline{B_r(o)}\subset E$.

On the other hand, suppose that $x\in E\setminus B_{2r}(o)$. For $r/3<s<r/2$, the ball
$B_s(x)$ lies outside $B_{3r/2}(o)$, where $h=-3n/r$. In particular,
$E\setminus B_s(x)$ is an admissible competitor for almost every such $s$.
Applying Lemma~\ref{lem:ball-perimeter-comparison} once more, we obtain
$$
 0\leq\mathcal J(E\setminus B_s(x))-\mathcal J(E)
 \leq\left(\frac ns-\frac{3n}{r}\right)\Vol(E\cap B_s(x)).
$$
The coefficient is negative, whereas $\Vol(E\cap B_s(x))>0$ since
$E$ is open and contains $x$. This contradiction proves $E\subset B_{2r}(o)$.
\end{proof}

\begin{proof}[Proof of Lemma~\ref{lem:prescribed-boundaries}]
According to Lemma~\ref{lem:direct-existence} and Lemma~\ref{lem:direct-boundary-properties},
there is an open minimizer $E$ with
$\overline{B_r(o)}\subset E\subset B_{2r}(o)$.
Its boundary lies in $M\setminus\overline{B_r(o)}$, where
$E$ is a local quasiminimizer of perimeter. The density estimates for
quasiminimizers \cite[Theorem~4.2]{KinnunenKorteLorentShanmugalingam}
give $$\Gamma=\partial E=\spt|D\chi_E|.$$
If $E\Delta F$ is compactly contained in a coordinate ball of radius $s$,
the volume term in \eqref{eq:exhaustion-almost-minimality} has the order $O(s^n)$.
Since the metric coefficients are smooth, \cite[Theorem~1.1]{Simmons}
gives $C^{1,\alpha}$ regularity of
$\Gamma$ away from a closed singular set of Hausdorff dimension at most $n-8$.
This set is empty for $n\leq7$.

The first variation of $\mathcal J$ gives the generalized mean
curvature vector of $|\Gamma|$ as $-h\nu_E$ almost everywhere, where
$\nu_E$ is the outward unit normal. Since $h$ is smooth on $M$,
elliptic regularity for the prescribed mean curvature equation implies
that the regular part of $\Gamma$ is smooth. The bound $|h|\leq3n/r$
gives \eqref{eq:direct-mean-curvature}, and we complete the proof.
\end{proof}

\begin{proof}[Proof of Proposition~\ref{prop:caccioppoli-exhaustion}]
Fix $o\in M$. For each $i\geq1$, Lemma~\ref{lem:prescribed-boundaries}
with $r=2^i$ gives an open Caccioppoli set $E_i$ such that
$$
 \overline{B_{2^i}(o)}\subset E_i\subset B_{2^{i+1}}(o).
$$
Thus we have $\overline{E_i}\subset E_{i+1}$ and that
 $\overline{E_i}$ is compact by the Hopf--Rinow theorem.
Since the balls $B_{2^i}(o)$ exhaust $M$, we have
$\bigcup_{i=1}^\infty E_i=M$.
Lemma~\ref{lem:prescribed-boundaries} also gives the desired regularity
of $\Gamma_i=\partial E_i=\spt|D\chi_{E_i}|$. By
\eqref{eq:direct-mean-curvature}, we have
$$
 \|\mathbf H_i\|_{L^\infty(\Gamma_i)}\leq\frac{3n}{2^i}\to0.
$$
We complete the proof.
\end{proof}

\section{An intrinsic spectral volume comparison theorem}
\label{sec:spectral-comparison}

In this section, we will extend the spectral volume comparison of Antonelli and Xu
\cite{AntonelliXu} to compact AM--PI spaces (see Definition~\ref{def:ampi}).
Fix $m\geq3$. Let $(X,\mathbf d,\boldsymbol{\mu})$ be a compact $m$-dimensional AM--PI space
with decomposition 
$$X=\mathcal R\sqcup\mathcal S$$ and connected regular
part $\mathcal R$. By definition, $\mathcal S$ is closed
and $\boldsymbol{\mu}(\mathcal S)=0$. The length metric $\mathbf d_\ell$
is defined to be the infimum of the $\mathbf d$-lengths of curves in $\mathcal R$
joining two points. It induces the subspace topology on $\mathcal R$.
Fix a smooth Riemannian manifold $(N,g)$ and a measure-preserving isometry
$$
 \iota\colon(\mathcal R,\mathbf d_\ell,\boldsymbol{\mu}|_{\mathcal R})
 \xrightarrow{\;\cong\;}(N,d_g,\mathrm{dvol}_g).
$$

In the following, we identify
$(\mathcal R,\mathbf d_\ell,\boldsymbol{\mu}|_{\mathcal R})$ with
the Riemannian manifold $(N,g)$ via $\iota$. We use $\mathcal H^s$
without a subscript to denote the Hausdorff measure for the ambient distance
$\mathbf d$. For $E\subset\mathcal R$, we use
$\mathcal H_g^{s}(E)$ to denote its Hausdorff measure with respect to
$\mathbf d_\ell$. We also denote $B_r(x)$ and $B_r^\ell(x)$ for balls in
$(X,\mathbf d)$ and $(\mathcal R,\mathbf d_\ell)$, respectively.

Throughout this section, we assume that $\mathcal H^{m-3}$ is $\sigma$-finite on
$\mathcal S$. Namely, $\mathcal S$ is a countable union of Borel sets
of finite $\mathcal H^{m-3}$-measure.

\begin{theorem}\label{thm:singularbishop}
Suppose that we have 
$$\Ric_g\geq\vartheta g$$ on $\mathcal R$ for a locally
Lipschitz function $\vartheta$ bounded above. If there is $\Lambda>0$ such that
$$
 \int_{\mathcal R}(|\nabla\phi|^2+\vartheta\phi^2)\,d\boldsymbol{\mu}
 \geq(m-1)\Lambda\int_{\mathcal R}\phi^2\,d\boldsymbol{\mu}
$$
for every $\phi\in C_c^\infty(\mathcal R)$, then
\begin{equation}\label{eq:singulararea}
 \boldsymbol{\mu}(X)\leq\Lambda^{-m/2}\omega_m.
\end{equation}
\end{theorem}

In the following, we work with the hypotheses of Theorem~\ref{thm:singularbishop}.
We first construct a positive supersolution to define the weighted
isoperimetric profile.

\begin{lemma}\label{lem:profile-weight}
For every $0<\kappa<\Lambda$ and every $\gamma>1$, there exists
$u\in W^{1,2}(X)\cap C^2(\mathcal R)$ such that
\begin{equation}\label{eq:profile-supersolution}
 -\gamma\Delta_{\mathcal R}u+\vartheta u\geq(m-1)\kappa u
 \,\,\text{on }\mathcal R\,\,\text{with}\,\,\inf_{\mathcal R}u=1.
\end{equation}
\end{lemma}

\begin{proof}
Since $\mathcal H^{m-3}$ is $\sigma$-finite on $\mathcal S$, we have
$$\mathcal H^{m-2}(\mathcal S)=0.$$ The spectral inequality gives
\eqref{eq:ampi-spectral-bound}
for $V=\vartheta-(m-1)\kappa$ with $\sigma=(m-1)(\Lambda-\kappa)>0$.
Lemma~\ref{lem:ampi-weight} therefore gives a function $u$ satisfying
\eqref{eq:profile-supersolution} with $\essinf_Xu=1$.
Since $\boldsymbol{\mu}(\mathcal S)=0$, we have $\essinf_{\mathcal R}u=1$.
Every nonempty open subset of $\mathcal R$ has positive Riemannian volume,
so the continuity of $u$ gives $\inf_{\mathcal R}u=1$.
\end{proof}

Fix
$$
 0<\kappa<\Lambda
 \,\,\text{and}\,\,
 1<\gamma<\frac{m-1}{m-2}.
$$
Choose $u$ as in Lemma~\ref{lem:profile-weight} and set
$$
 \alpha=\frac{2\gamma}{m-1}.
$$
For a set $F$ of locally finite perimeter in $\mathcal R$, define
$$
 V_u(F)=\int_F u^\alpha\,d\boldsymbol{\mu}\,\,\text{and}
 \,\,
 \mathcal P_u(F)=\int_{\partial^*F}u^\gamma\,d\mathcal H_g^{m-1},
$$
where $\partial^*F$ is the reduced boundary of $F$ in $\mathcal R$.
Since we have $u\in L^2(X)$ and $\alpha<2$, the total weighted volume
$V_u(\mathcal R)$ is finite. In the following, we consider the weighted isoperimetric profile
$$
 I(v)=\inf\{\mathcal P_u(F):V_u(F)=v\}\,\,\text{with}\,\, 0\leq v\leq V_u(\mathcal R).
$$

\begin{lemma}\label{lem:profile-regions}
For every $v\in\bigl(0,V_u(\mathcal R)\bigr)$, the infimum
$I(v)$ is attained. The weighted isoperimetric profile $I$ is positive and continuous on 
$\bigl(0,V_u(\mathcal R)\bigr)$. If $F$ is any minimizer attaining $I(v)$, there are constants
$C,r_0>0$, depending on $F$, such that
\begin{equation}\label{eq:profile-local-comparison}
 \mathcal P_u(F)\leq\mathcal P_u(E)+C|V_u(E)-v|
\end{equation}
whenever $E$ has locally finite perimeter in $\mathcal R$ and
$F\Delta E\subset B_r(x)$ for some $x\in X$ and $0<r<r_0$.
\end{lemma}

\begin{proof}
Fix $0<v<V_u(\mathcal R)$. Choose finitely many relatively compact
normal balls $B_{r_i}^\ell(x_i)$ whose union has weighted volume greater
than $v$. Since the geodesic spheres have zero volume, the weighted volume
of $\bigcup_i B_{tr_i}^\ell(x_i)$ varies continuously from zero to a
value greater than $v$ as $t$ increases from zero to one. The union
of volume $v$ has finite weighted perimeter, so $I(v)<\infty$.

Let $\{F_j\}_{j=1}^\infty$ be a minimizing sequence for $I(v)$ with
$V_u(F_j)=v$ for every $j$. Since we have $u\geq1$, the unweighted perimeters
of $F_j$ are bounded. By the compactness theorem for functions of
bounded variation on an exhaustion of $\mathcal R$, a subsequence of
$\{\chi_{F_j}\}_{j=1}^\infty$ converges almost everywhere to $\chi_F$
for a set $F$ of locally finite perimeter. Since $u^\alpha\in L^1(X)$,
the dominated convergence theorem gives $V_u(F)=v$. Lower semicontinuity
of weighted perimeter gives $\mathcal P_u(F)\leq I(v)$, so $F$ attains $I(v)$.

Since $\mathcal R$ is connected, a subset of zero perimeter has zero or
full volume. As $0<v<V_u(\mathcal R)$, every minimizer attaining $I(v)$
has positive perimeter. Hence we have $I(v)>0$, verifying the positivity.

For continuity, we verify that $I$ is both upper and lower semicontinuity.
Let $F$ be any minimizer at volume $v$, and let $B\Subset\mathcal R$ be
a coordinate ball with $\mathcal H_g^{m-1}(B\cap\partial^*F)>0$.
Since $u^\alpha>0$ and the outward unit normal $\nu_F$ has length one,
the vector measure 
$$u^\alpha\nu_F\,\mathcal H_g^{m-1}\llcorner(B\cap\partial^*F)$$
is nonzero. Choose a smooth vector field $Y$ compactly supported in $B$
with nonzero pairing against the above measure. For its flow $\Phi_t$, the
Gauss--Green formula gives
$$
 \left.\frac{d}{dt}V_u(\Phi_t(F))\right|_{t=0}
 =\int_{B\cap\partial^*F}u^\alpha\langle Y,\nu_F\rangle_g\,
 d\mathcal H_g^{m-1}\ne0.
$$
The area formula gives differentiability of $\mathcal P_u(\Phi_t(F))$.
By the inverse function theorem, there are $C,\varepsilon>0$ such that,
for every $|\delta|<\varepsilon$, one can choose $t$ with
\begin{equation}\label{eq:profile-volume-adjustment}
 V_u(\Phi_t(F))=v+\delta
 \,\,\text{and}\,\,
 \mathcal P_u(\Phi_t(F))\leq I(v)+C|\delta|.
\end{equation}
Hence $I(v+\delta)\leq I(v)+C|\delta|$, which gives the upper semicontinuity at $v$.

Suppose that $I$ is not lower semicontinuous at $v$. Then there are
$\eta>0$ and a sequence $\{v_j\}_{j=1}^\infty$ converging to $v$ such
that $I(v_j)\leq I(v)-\eta$ for every $j$. At each volume $v_j$,
choose a minimizer $E_j$. Since $\mathcal P_u(E_j)=I(v_j)$, their
weighted perimeters are uniformly bounded. By the compactness argument
above, a subsequence of their characteristic functions converges almost
everywhere on $\mathcal R$ to $\chi_E$ for a set $E$ of locally finite
perimeter. The dominated convergence theorem gives $V_u(E)=v$, and lower
semicontinuity of weighted perimeter gives $\mathcal P_u(E)\leq I(v)-\eta$.
This contradicts the definition of $I(v)$. Thus $I$ is also lower
semicontinuous at $v$, which proves continuity.

To prove \eqref{eq:profile-local-comparison}, we choose two coordinate balls
$U_1,U_2\Subset\mathcal R$ with disjoint closures, each centered at a
point of $\partial^*F$. In each $U_i$, we choose a compactly supported
vector field $Y_i$ as above. Choose $C,\varepsilon>0$ so that, for
$i=1,2$ and every $|\delta|<\varepsilon$, the flow of $Y_i$ satisfies
\eqref{eq:profile-volume-adjustment} at some time $t=t(i,\delta)$.
By the H\"older inequality and the upper Ahlfors bound, we have
$$
 V_u(B_r(x)\cap\mathcal R)
 \leq \|u\|_{L^2(X)}^\alpha\boldsymbol{\mu}(B_r(x))^{1-\alpha/2}
 \to0
$$
uniformly in $x\in X$ as $r\downarrow0$. Choose $r_0>0$ so that
$2r_0$ is less than the $\mathbf d$-distance between $\overline{U_1}$
and $\overline{U_2}$ and $V_u(B_{r_0}(x)\cap\mathcal R)<\varepsilon$
for every $x\in X$.

Let $E$ be a competitor with $F\Delta E\subset B_r(x)$, where $x\in X$
and $0<r<r_0$. We may assume $\mathcal P_u(E)<\infty$. Since $B_r(x)$
has diameter at most $2r$, it misses one of the balls $U_i$. Thus $E$
and $F$ agree on $U_i$. Set $\delta=v-V_u(E)$. Since
$|\delta|\leq V_u(B_r(x)\cap\mathcal R)<\varepsilon$, we may apply
\eqref{eq:profile-volume-adjustment} to the flow $\Phi_t$ of $Y_i$ with
this value of $\delta$.
The flow acts only in $U_i$, so it changes the weighted volume and
perimeter of $E$ by the same amounts as for $F$. Hence we have
$$
 V_u(\Phi_t(E))=v
 \,\,\text{and}\,\,
 \mathcal P_u(\Phi_t(E))\leq\mathcal P_u(E)+C|\delta|.
$$
Comparing $F$ with $\Phi_t(E)$ gives \eqref{eq:profile-local-comparison}.
\end{proof}

Fix a minimizer $F$ attaining $I(v)$ with $0<v<V_u(\mathcal R)$. The local
perimeter comparison gives the following regularity for the boundary of $F$.

\begin{lemma}\label{lem:profile-regularity}
Let $$\Gamma=\spt_{\mathcal R}|D\chi_F|.$$ 
Then there is a relatively closed
set $Z_F\subset\Gamma$ such that $\Sigma=\Gamma\setminus Z_F$ is a
$C^2$ hypersurface. The set $Z_F$ is empty for $m\leq7$ and has Hausdorff
dimension at most $m-8$ with respect to $\mathbf d_\ell$ for $m\geq8$.
\end{lemma}

\begin{proof}
On a relatively compact coordinate chart, the error in
\eqref{eq:profile-local-comparison} has the order $O(r^m)$ whenever $F\Delta E$
lies in a coordinate ball of radius $r$.
The weighted perimeter is the perimeter for the metric
$u^{2\gamma/(m-1)}g$. By \cite[Theorem~1.1]{Simmons},
$\Gamma$ is $C^{1,\theta}$ for some
$\theta>0$ outside a relatively closed set $Z_F$ with the stated
dimension bound, and is regular everywhere when $m\leq7$. On the regular graphs the first
variation under the volume constraint gives
$$
 \mathbf H=c\,u^{\alpha-\gamma}-\gamma\,\partial_{\nu_F}\log u,
$$
where $c$ is the Lagrange multiplier for the constraint $V_u(F)=v$
and $\nu_F$ is the outward unit normal.
Interior Schauder regularity for this equation gives $C^2$ regularity
of $\Sigma$.
\end{proof}

To estimate the profile $I$, we use the second variation of
$\mathcal P_u-cV_u$ at $F$. For a variation with outward normal speed
$f\in C_c^1(\Sigma)$, this second variation contains the term
$$
 \int_\Sigma u^\gamma|\nabla_\Sigma f|^2\,d\mathcal H_g^{m-1}.
$$
For our purpose, we will take $f=u^{-\gamma}\chi$, where $\chi$ is a cutoff function
vanishing near $K=\mathcal S\cup Z_F$. The term quadratic in
$\nabla_\Sigma\chi$ is then
$$
 \int_\Sigma u^{-\gamma}|\nabla_\Sigma\chi|^2\,d\mathcal H_g^{m-1}.
$$
To choose $\chi$ with arbitrarily small cutoff error, we first estimate
the measure 
$$\nu=u^{-\gamma}\mathcal H_g^{m-1}|_\Sigma$$ on balls
centered at points of $K$.

\begin{lemma}\label{lem:profile-growth}
There are constants $C,r_0>0$ such that, for every $x\in K$ and
$0<r<r_0$, we have
$$
 \nu(B_r(x))\leq Cr^{m-1}.
$$
\end{lemma}

\begin{proof}
We first derive a weighted boundary estimate from the local perimeter
inequality by removing superlevel sets from $F$.
By \eqref{eq:profile-local-comparison} and the upper Ahlfors bound,
we can choose $C>0$ and $0<r_0\leq1$ such that, for every $x\in K$
and $0<r<r_0$,
$$
 \mathcal P_u(F)\leq\mathcal P_u(E)+C|V_u(E)-v|
$$
for every set $E$ of locally finite perimeter in $\mathcal R$ with
$F\Delta E\subset B_{2r}(x)$, and
$$
 \boldsymbol{\mu}(B_{2r}(x))\leq Cr^m.
$$

Fix $x\in K$ and $0<r<r_0$.
We claim that every locally Lipschitz function $\psi:\mathcal R\to[0,1]$
vanishing outside $B_{2r}(x)$ satisfies
\begin{equation}\label{eq:profile-weighted-boundary}
 \int_\Sigma\psi u^\gamma\,d\mathcal H_g^{m-1}
 \leq\int_{\mathcal R}u^\gamma|\nabla\psi|\,d\boldsymbol{\mu}
      +C\int_F\psi u^\alpha\,d\boldsymbol{\mu}.
\end{equation}
We may assume that the right-hand side is finite.
Set $$E_t=F\setminus\{\psi>t\}$$ for $0<t<1$.
For almost every $t$, the coarea formula gives locally finite perimeter
of $\{\psi>t\}$, and $\mathcal H_g^{m-1}(\Sigma\cap\{\psi=t\})=0$
because $\mathcal H_g^{m-1}(\Sigma)<\infty$.
The perimeter formula for a set difference \cite[Theorem~16.3]{Maggi},
applied in local coordinates on $\mathcal R$, gives
$$
 \mathcal P_u(E_t)-\mathcal P_u(F)
 \leq\int_{\{\psi=t\}}u^\gamma\,d\mathcal H_g^{m-1}
     -\int_{\Sigma\cap\{\psi>t\}}u^\gamma\,d\mathcal H_g^{m-1}.
$$
Since $F\Delta E_t\subset B_{2r}(x)$ and
$v-V_u(E_t)=V_u(F\setminus E_t)$, \eqref{eq:profile-local-comparison}
applies to $E_t$. We add $CV_u(F\setminus E_t)$ to the perimeter
estimate and integrate over $0<t<1$. By the coarea formula on
$\mathcal R$ and Tonelli's theorem on $\Sigma$ and $F$, we obtain
$$
\begin{aligned}
 0&\leq\int_0^1\bigl(\mathcal P_u(E_t)-\mathcal P_u(F)
       +CV_u(F\setminus E_t)\bigr)\,dt\\
 &\leq\int_{\mathcal R}u^\gamma|\nabla\psi|\,d\boldsymbol{\mu}
       +C\int_F\psi u^\alpha\,d\boldsymbol{\mu}\\
 &\quad-\int_\Sigma\psi u^\gamma\,d\mathcal H_g^{m-1}.
\end{aligned}
$$
This proves \eqref{eq:profile-weighted-boundary}.

To recover the measure $\nu$, let $\eta$ be the radial cutoff function
equal to one on $B_r(x)$, zero outside $B_{2r}(x)$, and linear in
$\mathbf d(x,\cdot)$ between the two radii. We take $
 \psi=\eta u^{-2\gamma}
$
so that
$\psi u^\gamma=\eta u^{-\gamma}$ in \eqref{eq:profile-weighted-boundary}.
Since $u\geq1$ and $\alpha\leq\gamma$, we have
$$
 u^\gamma|\nabla\psi|\leq|\nabla\eta|+2\gamma\eta|\nabla\log u|
 \,\,\text{and}\,\,\psi u^\alpha\leq\eta.
$$
Using the gradient bound $|\nabla\eta|\leq r^{-1}$ and the upper Ahlfors bound in
\eqref{eq:profile-weighted-boundary}, we obtain
\begin{equation}\label{eq:profile-growth-reduction}
 \nu(B_r(x))\leq Cr^{m-1}
       +2\gamma\int_{\mathcal R}\eta|\nabla\log u|\,d\boldsymbol{\mu}.
\end{equation}
It remains to bound the integral on the right by $Cr^{m-1}$.

Choose $a\geq0$ with $\gamma a\geq\sup_{\mathcal R}\vartheta-(m-1)\kappa$.
Then \eqref{eq:profile-supersolution} gives
$$-\Delta_{\mathcal R}u+au\geq0$$ on $\mathcal R$.
Since we have $u\geq1$ and $1\in W^{1,2}(X)$ by Lemma~\ref{lem:ampi-constant},
the Sobolev chain rule and multiplication by $\eta^2$ give
$\eta^2/u\in W^{1,2}(X)$. Approximation by nonnegative functions in
$C_c^\infty(\mathcal R)$ allows us to test the inequality with $\eta^2/u$.
Young's inequality gives
$$
 \int_{\mathcal R}\eta^2|\nabla\log u|^2\,d\boldsymbol{\mu}
 \leq\int_{\mathcal R}(4|\nabla\eta|^2+2a\eta^2)\,d\boldsymbol{\mu}
 \leq Cr^{m-2}.
$$
The Cauchy--Schwarz inequality now yields
$$
 \int_{\mathcal R}\eta|\nabla\log u|\,d\boldsymbol{\mu}
 \leq\boldsymbol{\mu}(B_{2r}(x))^{1/2}(Cr^{m-2})^{1/2}
 \leq Cr^{m-1}.
$$
Substitution in \eqref{eq:profile-growth-reduction} proves the desired bound.
\end{proof}

The estimate in Lemma~\ref{lem:profile-growth} allows us to construct
cutoff functions near $K$ with arbitrarily small Dirichlet energy with 
respect to the measure $\nu$. In the following,
 $\operatorname{lip}$ denotes the pointwise upper Lipschitz constant
with respect to $\mathbf d$.

\begin{lemma}\label{lem:metric-hausdorff-cutoffs}
Let $(X,\mathbf d)$ be a compact metric space with a finite Borel measure
$\nu$. Let $K\subset X$ be compact and assume that $\mathcal H^{m-3}$ is
$\sigma$-finite on $K$, where $m\geq3$. Suppose that for some $C,r_0>0$,
$$
 \nu(B_r(y))\leq Cr^{m-1}\,\,\text{for every }y\in K
 \,\,\text{and }0<r<r_0.
$$
For every open neighborhood $U$ of $K$ and every $\varepsilon>0$,
there is a Lipschitz function $\chi:X\to[0,1]$, zero near $K$ and
equal to one on $X\setminus U$, such that
$$
 \int_X(\operatorname{lip}\chi)^2\,d\nu<\varepsilon.
$$
\end{lemma}

\begin{proof}
If $K$ is empty, take $\chi=1$. For $y\in K$ and $0<r<r_0/3$,
let $\eta_{y,r}$ equal one on $B_r(y)$, vanish outside $B_{2r}(y)$,
and decrease linearly in the distance to $y$ between the two radii. Then
$$
 \int_X(\operatorname{lip}\eta_{y,r})^2\,d\nu
 \leq r^{-2}\nu(B_{3r}(y))\leq 3^{m-1}Cr^{m-3}.
$$
The energy of a finite maximum of Lipschitz functions is at most the sum
of their energies, since its pointwise slope is at most the largest of
their slopes.

First let $L\subset K$ be compact with $\mathcal H^{m-3}(L)<\infty$.
By the definition of Hausdorff measure and compactness, $L$ has a finite
cover by balls $B_{r_i}(y_i)$ centered on $L$, with arbitrarily small
radii and
$$
 \sum_i r_i^{m-3}\leq C_m\bigl(1+\mathcal H^{m-3}(L)\bigr).
$$
The maximum of the corresponding cutoff functions equals one near $L$ and has
energy at most $C_L$, independently of the radii. Its support can thus
be taken in any prescribed neighborhood of $L$. Choose $N$ such cutoff functions
$\eta_1,\ldots,\eta_N$, the first supported in $U$ and each subsequent
one supported where its predecessor is locally equal to one. Their average
$\eta$ equals one near $L$. At each point at most one cutoff function has nonzero
slope, so
$$
 \int_X(\operatorname{lip}\eta)^2\,d\nu
 \leq \frac{1}{N^2}\sum_{i=1}^N\int_X(\operatorname{lip}\eta_i)^2\,d\nu
 \leq \frac{C_L}{N}.
$$

By $\sigma$-finiteness and inner regularity on each finite-measure
Borel piece, we may write
$$
 K=Z\cup\bigcup_{i=1}^{\infty}L_i,
$$
where each $L_i$ is compact with finite $\mathcal H^{m-3}$-measure and
$\mathcal H^{m-3}(Z)=0$. When $m=3$, each $L_i$ is finite and $Z$ is empty.
For each $L_i$, choose a cutoff function as above
supported in $U$ with energy less than $2^{-i-1}\varepsilon$.
Cover $Z$ by countably many balls $B_r(y)$ centered on $Z$, with $3r<r_0$ and
$\overline{B_{2r}(y)}\subset U$. Since $\mathcal H^{m-3}(Z)=0$, we can
make $\sum r^{m-3}$ so small that the cutoff functions for these balls have
total energy less than $\varepsilon/2$.

The open sets where these cutoff functions are locally equal to one cover $K$.
Choose a finite subcover and let $h$ be the maximum of the corresponding
cutoff functions. Then $h$ equals one near $K$, is supported in $U$, and has
energy less than $\varepsilon$.
Taking $\chi=1-h$ proves the lemma.
\end{proof}

We apply Lemma~\ref{lem:metric-hausdorff-cutoffs} to make the error in the
second variation on $\Sigma$ tend to zero. This gives the following
differential inequality for the weighted isoperimetric profile $I$.

\begin{lemma}\label{lem:profile-inequality}
Fix $0<v<V_u(\mathcal R)$. If $\psi$ is $C^2$ near $v$,
$\psi(v)=I(v)$, and $\psi\leq I$ near $v$, then we have
$$
 I(v)\psi''(v)\leq-\frac{\psi'(v)^2}{m-1}-(m-1)\kappa.
$$
\end{lemma}

\begin{proof}
Choose a minimizer $F$ attaining $I(v)$ and let $\Sigma$ be its regular boundary.
For any $\chi\in C_c^1(\Sigma)$, vary $F$ with outward normal speed
$u^{-\gamma}\chi$. Denote the weighted perimeter and volume of the
varied set by $P(t)$ and $V(t)$, respectively. Since
$$P(t)\geq I(V(t))\geq\psi(V(t))$$ for small $t$, with equality at
$t=0$, we have 
$$P'(0)=\psi'(v)V'(0)$$ and
$$
 \psi''(v)V'(0)^2\leq P''(0)-\psi'(v)V''(0).
$$
Denote
$$
 q(\chi):=V'(0)=\int_\Sigma u^{\alpha-\gamma}\chi\,d\mathcal H_g^{m-1}.
$$
Using $\Ric_g\geq\vartheta g$, \eqref{eq:profile-supersolution},
$\alpha=2\gamma/(m-1)$, and $\gamma<(m-1)/(m-2)$, the second variation
calculation in
\cite[Appendix~A, equations~(84)--(86)]{AntonelliXu} gives
$$
\begin{aligned}
 P''(0)-\psi'(v)V''(0)
 &\leq\int_\Sigma|\nabla_\Sigma\chi|^2\,d\nu\\
 &\quad-\left(\frac{\psi'(v)^2}{m-1}+(m-1)\kappa\right)
 \int_\Sigma u^{2\alpha-3\gamma}\chi^2\,d\mathcal H_g^{m-1}.
\end{aligned}
$$
Since $I(v)=\int_\Sigma u^\gamma\,d\mathcal H_g^{m-1}$, the
Cauchy--Schwarz inequality gives
$$q(\chi)^2\leq I(v)\int_\Sigma u^{2\alpha-3\gamma}\chi^2\,d\mathcal H_g^{m-1}.$$
Combining these estimates yields
\begin{equation}\label{eq:profile-cutoff-variation}
 \left(I(v)\psi''(v)+\frac{\psi'(v)^2}{m-1}+(m-1)\kappa\right)q(\chi)^2
 \leq I(v)\int_\Sigma|\nabla_\Sigma\chi|^2\,d\nu.
\end{equation}
By $W^{1,2}$ approximation on compact subsets of $\Sigma$, this inequality
also holds for compactly supported Lipschitz functions $\chi$ on $\Sigma$.

Since $Z_F$ is relatively closed in $\mathcal R$, the set
$K=\mathcal S\cup Z_F$ is closed in $X$ and hence compact.
Moreover, $u\geq1$ gives $\nu(X)\leq I(v)<\infty$.
Since $I(v)>0$, we can fix a nonempty open set $D\Subset\Sigma$.
Choose an open neighborhood $U$ of $K$ disjoint from $\overline D$.
By Lemma~\ref{lem:profile-regularity} and
$\mathbf d\leq\mathbf d_\ell$, we have $\mathcal H^{m-3}(Z_F)=0$.
Thus $\mathcal H^{m-3}$ is $\sigma$-finite on $K$.
Given $\varepsilon>0$, Lemmas~\ref{lem:profile-growth}
and~\ref{lem:metric-hausdorff-cutoffs} give a Lipschitz cutoff function
$0\leq\chi\leq1$ on $X$, zero near $K$ and equal to one outside $U$,
with
\begin{equation}\label{eq:profile-cutoff-energy}
 \int_X(\operatorname{lip}\chi)^2\,d\nu<\varepsilon.
\end{equation}
Its restriction to $\Sigma$ has compact support and satisfies
$|\nabla_\Sigma\chi|\leq\operatorname{lip}\chi$ almost everywhere.
Since $\chi=1$ on $D$, we have
$$
 q(\chi)\geq\int_D u^{\alpha-\gamma}\,d\mathcal H_g^{m-1}>0.
$$
The lower bound is independent of $\varepsilon$.
Substituting \eqref{eq:profile-cutoff-energy} into
\eqref{eq:profile-cutoff-variation} and letting $\varepsilon\to0$
gives the desired inequality.
\end{proof}

\begin{proof}[Proof of Theorem~\ref{thm:singularbishop}]
Fix $x\in\mathcal R$. In normal coordinates at $x$, the continuity of $u$
gives,
$$
\begin{aligned}
 V_u(B_s^\ell(x))&=\frac{\omega_{m-1}}{m}u(x)^\alpha s^m(1+o(1))\\
 \text{and}\,\,\mathcal P_u(B_s^\ell(x))
 &=\omega_{m-1}u(x)^\gamma s^{m-1}(1+o(1))\,\,\text{as}\,\,s\downarrow 0.
\end{aligned}
$$
For each sufficiently small $v>0$, we can choose $s$ so that $V_u(B_s^\ell(x))=v$.
By the definition of $I$, we have 
$$I(v)\leq\mathcal P_u(B_s^\ell(x)).$$
Dividing by $v^{(m-1)/m}$ and using $\alpha=2\gamma/(m-1)$ gives
$$
\begin{aligned}
 \limsup_{v\downarrow0}\frac{I(v)}{v^{(m-1)/m}}
 &\leq\lim_{s\downarrow0}
 \frac{\mathcal P_u(B_s^\ell(x))}{V_u(B_s^\ell(x))^{(m-1)/m}}\\
 &=\bigl(m^{m-1}\omega_{m-1}\bigr)^{1/m}u(x)^{\gamma(m-2)/m}.
\end{aligned}
$$
Since $\inf_{\mathcal R}u=1$, taking the infimum over $x$ yields
\begin{equation}\label{eq:profile-at-zero}
 \limsup_{v\downarrow0}\frac{I(v)}{v^{(m-1)/m}}
 \leq\bigl(m^{m-1}\omega_{m-1}\bigr)^{1/m}.
\end{equation}
In particular, we see $I(v)\to0=I(0)$ as $v\downarrow0$.

By Lemma~\ref{lem:profile-regions}, $I$ is positive and continuous on
$(0,V_u(\mathcal R))$. Applying the one-dimensional comparison in
\cite[Lemma~4]{AntonelliXu} to $I$, with the
differential inequality of Lemma~\ref{lem:profile-inequality} on
$(0,V_u(\mathcal R))$ and the estimate \eqref{eq:profile-at-zero} at $v=0$, gives
$$V_u(\mathcal R)\leq\kappa^{-m/2}\omega_m.$$ Since $u\geq1$ and
$\boldsymbol{\mu}(\mathcal S)=0$, we have
$$
 \boldsymbol{\mu}(X)=\boldsymbol{\mu}(\mathcal R)
 \leq\int_{\mathcal R}u^\alpha\,d\boldsymbol{\mu}
 =V_u(\mathcal R).
$$
Letting $\kappa\uparrow\Lambda$ yields \eqref{eq:singulararea}.
\end{proof}

\section{Proof of the main theorem}
\label{sec:volume}

We first use radial projection to bound the area of a portion of a
geodesic sphere in terms of the area of a countably rectifiable set.
For the case of smooth separating hypersurfaces, see
\cite{WeiXuZhang,ZhouZhu} for instance.

\begin{lemma}\label{lem:radial}
Let $(M^{m+1},g)$ be a connected complete Riemannian manifold with
$\Ric_M\geq0$, and fix $q\in M$. Let
$\Sigma\subset M\setminus\{q\}$ be a nonempty compact, countably
$m$-rectifiable set of finite $\mathcal H_g^m$-measure, and set
$\mathfrak{d}=d(q,\Sigma)>0$. Suppose that $A\subset M$ is measurable and that,
for almost every $x\in A\setminus\operatorname{Cut}(q)$, the
minimizing segment from $q$ to $x$ meets $\Sigma$. Then for almost
every $s>0$,
\begin{equation}\label{eq:abstractradial}
 \mathcal H_g^m(A\cap\partial B_s(q))
 \leq \left(\frac{s}{\mathfrak{d}}\right)^m\mathcal H_g^m(\Sigma).
\end{equation}
\end{lemma}

\begin{proof}
Let $S_qM$ be the unit sphere in $T_qM$ with its induced round metric
and area measure $d\theta$. For $\theta\in S_qM$, let $c(\theta)$
denote the cut time. Write $J(r,\theta)$ for the polar volume density,
so that $dV_g=J(r,\theta)\,dr\,d\theta$ for $0<r<c(\theta)$.
For $s>0$, set
$$
 \Theta_s=\{\theta\in S_qM:s<c(\theta),\,\,\exp_q(s\theta)\in A\}.
$$

The cut locus has zero volume. By hypothesis, the same holds for
the set of points in $A\setminus\operatorname{Cut}(q)$ whose minimizing
geodesic segment from $q$ is disjoint from $\Sigma$.
The coarea formula for $d(q,\cdot)$ implies that both sets intersect
$\partial B_s(q)$ in sets of zero $\mathcal H_g^m$-measure for almost
every $s>0$. 

Therefore, for each such $s$, polar coordinates give
$$
 \mathcal H_g^m(A\cap\partial B_s(q))
 =\int_{\Theta_s}J(s,\theta)\,d\theta.
$$
Moreover, for $d\theta$-almost every $\theta\in\Theta_s$,
$$
 \Sigma\cap\{\exp_q(r\theta):\mathfrak{d}\leq r\leq s\}\neq\varnothing.
$$
Fix such an $s$ and define
$$
 \Sigma_s=\Sigma\cap\{\exp_q(r\theta):
          \theta\in\Theta_s,\,\,\mathfrak{d}\leq r\leq s\}.
$$

Consider the smooth radial projection from
$M\setminus(\{q\}\cup\operatorname{Cut}(q))$ to $S_qM$, given by
$$
 \Pi_q(\exp_q(r\theta))=\theta.
$$
For $\mathcal H_g^m$-almost every $z=\exp_q(r\theta)\in\Sigma_s$,
its tangential Jacobian satisfies
$$
 J_m^\Sigma\Pi_q(z)
 =\frac{|\langle \nu(z),\partial_r\rangle|}{J(r,\theta)}
 \leq\frac{1}{J(r,\theta)},
$$
where $\nu(z)$ is a unit normal to the approximate tangent plane
of $\Sigma$ at $z$. Since $\Ric_M\geq0$, the Bishop comparison theorem
for the polar volume density gives $J(s,\theta)/J(r,\theta)\leq(s/r)^m$
for $0<r\leq s<c(\theta)$.
Since $\Sigma_s\cap\Pi_q^{-1}(\theta)\neq\varnothing$ for
$d\theta$-almost every $\theta\in\Theta_s$, the area formula for
$\Pi_q|_{\Sigma_s}$ yields
$$
\begin{aligned}
 \int_{\Theta_s}J(s,\theta)\,d\theta
 &\leq\int_{\Sigma_s}J(s,\Pi_q z)J_m^\Sigma\Pi_q(z)
                   \,d\mathcal H_g^m(z)\\
 &\leq (s/\mathfrak{d})^m\mathcal H_g^m(\Sigma).
\end{aligned}
$$
This proves \eqref{eq:abstractradial}.
\end{proof}

We need one component of the boundary that separates a given interior
point from a given exterior point.

\begin{lemma}\label{lem:separating-component}
Let $(M^{m+1},g)$ be a connected orientable Riemannian manifold with
$m\geq1$ and $H_m(M;\mathbb Z)=0$. Let $E\subset M$ be a precompact
open Caccioppoli set with $\Gamma=\partial E=\spt|D\chi_E|$. Suppose
that $\mathcal S\subset\Gamma$ is closed, $\mathcal H_g^{m-1}(\mathcal S)=0$,
and $\Gamma\setminus\mathcal S$ is a smooth embedded hypersurface.
For every $o\in E$ and $q\in M\setminus\overline E$, there is a connected
component $\mathcal R_0$ of $\Gamma\setminus\mathcal S$ such that $o$ and
$q$ lie in different components of $M\setminus\overline{\mathcal R_0}$.
\end{lemma}

\begin{proof}
Choose a tubular neighborhood $U$ of a smooth embedded arc from $o$ to $q$,
with its end sections contained in $E$ and $M\setminus\overline E$,
respectively. Let $\pi:U\to D^m$ be the projection onto the transverse
$m$-disk. Since $\mathcal H_g^{m-1}(\mathcal S)=0$, the set
$\pi(\mathcal S\cap U)$ has zero $m$-dimensional measure.
By Sard's theorem applied to $\pi|_{(\Gamma\setminus\mathcal S)\cap U}$,
almost every fiber of $\pi$ is transverse to $\Gamma\setminus\mathcal S$.
Choose a fiber that is transverse and disjoint from $\mathcal S$.
Connecting its endpoints to $o$ and $q$ within the end sections gives
a piecewise smooth path $\gamma:[0,1]\to M$ disjoint from $\mathcal S$ and transverse
to $\Gamma\setminus\mathcal S$. The set $\gamma^{-1}(\Gamma)$ is compact
and discrete, so $\gamma$ intersects $\Gamma$ at finitely many points.

The equality $\Gamma=\partial E=\spt|D\chi_E|$ implies that $E$ lies
locally on one side of $\Gamma\setminus\mathcal S$.
Orient $\Gamma\setminus\mathcal S$ by the outward unit normal of $E$.
Since $\gamma$ joins $o\in E$ to $q\in M\setminus\overline E$, its
signed intersection number with $\Gamma\setminus\mathcal S$ is one.
Consequently, there is a connected component $\mathcal R_0$ of
$\Gamma\setminus\mathcal S$ whose signed intersection number with
$\gamma$ is nonzero.

Set $T=(\partial[E])\llcorner\mathcal R_0$. This is an integer
rectifiable current of finite mass and compact support.
The component $\mathcal R_0$ is a relatively closed oriented hypersurface
without boundary in $M\setminus\mathcal S$. Thus the flat chain
$\partial T$ is supported in $\mathcal S$.
Since $\mathcal H_g^{m-1}(\mathcal S)=0$, 
we have $\partial T=0$ from \cite[4.1.20]{Federer}.

By $H_m(M;\mathbb Z)=0$ and \cite[4.4.1]{Federer}, the integral cycle
$T$ is null-homologous. Its algebraic intersection number with every
closed curve transverse to $\mathcal R_0$ and disjoint from
$\spt T\setminus\mathcal R_0$ therefore vanishes.
If $o$ and $q$ belonged to the same component of
$M\setminus\overline{\mathcal R_0}$, there would be a piecewise smooth
path from $q$ to $o$ in that component. Concatenating this path with
$\gamma$ would give a closed curve with the same nonzero intersection
number with $T$ as $\gamma$, a contradiction.
\end{proof}

\begin{proof}[Proof of Theorem~\ref{thm:main}]
By \cite[Theorems~1.5 and~1.9]{ZhouZhu}, it suffices to consider
$n\geq8$. In the following, we write $m=n-1$ and fix $p\in M$. The proof
 will be divided into two cases.

\medskip
\noindent\textbf{Case 1.} $M=N\times\mathbb R$.

For each $x\in N$, choose
$t$ with $(x,t)\notin K$. For every unit vector $e\in T_xN$, we have
$$
 \Ric_N(e,e)=\biRic_M(\partial_t,e)\geq m-1.
$$
By the Bonnet--Myers theorem and the Bishop--Gromov volume comparison theorem,
$N$ is compact with diameter $D\leq\pi$ and volume at most $\omega_m$.
If $p=(x_0,t_0)$, then $B_R(p)$ lies in $N\times(t_0-R,t_0+R)$.
For $R>D$ it contains $N\times(t_0-R+D,t_0+R-D)$. Consequently, we obtain
$$
 \Vol B_R(p)\leq2R\Vol(N)
 \,\,\text{and}\,\,
 \lim_{R\to\infty}\frac{\Vol B_R(p)}R=2\Vol(N).
$$
This proves both bounds since $e(M)=2$.

\medskip
\noindent\textbf{Case 2.} $M$ does not split as a product with a line.

By the splitting theorem \cite{CheegerGromoll}, $M$ has one
end. We first treat the orientable case, for which
$H_m(M;\mathbb Z)=0$ by \cite[Theorem~1.1]{ShenSormani}.
In this case, we construct a hypersurface of area at most $\omega_m+o(1)$
separating a large ball about $p$ from infinity.

Choose $r_1\geq6n$ so that $K\subset B_{r_1}(p)$, and put
$r_j=4^{j-1}r_1$ for $j\geq1$. By Lemma~\ref{lem:prescribed-boundaries},
there are precompact open Caccioppoli sets $\{E_j\}_{j=1}^\infty$ with
$$
 \overline{B_{r_j}(p)}\subset E_j\subset B_{2r_j}(p).
$$
In particular, $\overline E_j\subset E_{j+1}$. Choose a smooth
approximation $\rho$ of $d(p,\cdot)$ with $|\rho-d(p,\cdot)|\leq1$
and $|\nabla\rho|\leq2$ by \cite{GreeneWu}. To prescribe the mean
curvature function between $E_j$ and $E_{j+1}$, set
\begin{equation}\label{eq:annular-linear-prescription}
 h_j=\frac1j\left(3-\frac{\rho}{r_j}\right).
\end{equation}
On $E_j$ we have $\rho\leq2r_j+1$, and outside $E_{j+1}$ we have
$\rho\geq4r_j-1$. Since $r_j\geq6nj$, it follows that
\begin{equation}\label{eq:annular-obstacle-bounds}
 h_j>\frac{3n}{r_j}\,\,\text{on }E_j
 \,\,\text{and}\,\,
 h_j<-\frac{3n}{r_{j+1}}\,\,\text{on }M\setminus E_{j+1}.
\end{equation}

Minimize the functional
$$
 \mathcal J_j(F)=\mathcal P(F)-\int_F h_j\,dV
$$
among Caccioppoli sets $F$ satisfying $E_j\subset F\subset E_{j+1}$
up to null sets. We call such sets \emph{admissible}. Since $h_j$
is bounded on $\overline E_{j+1}$, a minimizing sequence has bounded
perimeter and lies in a fixed compact set. By the compactness theorem
for sets of finite perimeter, a subsequence converges in $L^1(M)$ to
an admissible set $\Omega_j$. Lower semicontinuity of perimeter and
continuity of the integral term give the minimality of $\Omega_j$.
Taking $(\Omega_j\cup E_j)\cap E_{j+1}$ as the representative of
$\Omega_j$, we have $E_j\subset\Omega_j\subset E_{j+1}$ pointwise.

We show that $\Omega_j$ also minimizes $\mathcal J_j$ under compactly
supported variations in $M\setminus\overline{B_{r_j}(p)}$. Let $F$ be a
Caccioppoli set with $F\Delta\Omega_j\Subset M\setminus\overline{B_{r_j}(p)}$
and set $G=F\cup E_j$. Then $E_j\setminus F$ is
compactly contained in $M\setminus\overline{B_{r_j}(p)}$. Applying
the submodularity inequality for perimeter to $F$ and $E_j$, and
\eqref{eq:exhaustion-almost-minimality} to $E_j$ and $F\cap E_j$, gives
$$
 \mathcal J_j(G)-\mathcal J_j(F)
 \leq\int_{E_j\setminus F}\left(\frac{3n}{r_j}-h_j\right)\,dV\leq0
$$
by \eqref{eq:annular-obstacle-bounds}. On the other hand,
$G\setminus E_{j+1}=F\setminus E_{j+1}$ is compactly contained in
$M\setminus\overline{B_{r_{j+1}}(p)}$, since $E_{j+1}$ contains
$\overline{B_{r_{j+1}}(p)}$. Applying the same perimeter inequalities
to $G$ and $E_{j+1}$, with $G\cup E_{j+1}$ as competitor for $E_{j+1}$,
gives
$$
 \mathcal J_j(G\cap E_{j+1})-\mathcal J_j(G)
 \leq\int_{G\setminus E_{j+1}}\left(h_j+\frac{3n}{r_{j+1}}\right)\,dV\leq0.
$$
Since $G\cap E_{j+1}$ is admissible, we conclude that
$\mathcal J_j(\Omega_j)\leq\mathcal J_j(F)$.

Since $h_j$ is locally bounded, the argument in the proof of
Lemma~\ref{lem:direct-boundary-properties} shows that $\Omega_j$ is
locally quasiminimizing for perimeter outside $\overline{B_{r_j}(p)}$.
Since $\Vol(E_j\setminus\Omega_j)=0$, the density argument in that proof
allows us to take $\Omega_j$ open, with
$\Gamma_j=\partial\Omega_j=\spt|D\chi_{\Omega_j}|$.
The inclusions $E_j\subset\Omega_j\subset\overline E_{j+1}$ now give
$$
 \Gamma_j\subset\overline{B_{8r_j}(p)}\setminus\overline{B_{r_j}(p)}.
$$
By \eqref{eq:annular-linear-prescription}, we have
\begin{equation}\label{eq:annular-prescription-bounds}
 \sup_{\Gamma_j}|h_j|\leq\frac6j
 \,\,\text{and}\,\,
 \sup_M|\nabla h_j|\leq\frac{2}{jr_j}\leq\frac1{j^2}.
\end{equation}

The local minimality of $\Omega_j$ gives \eqref{eq:ampi-local-minimality}
with $E=\Omega_j$ and $h=h_j$ in a neighborhood of $\Gamma_j$.
By Proposition~\ref{prop:ampi-pmc}, $\Gamma_j$ is a compact AM--PI space
and 
$$\dim_{\mathcal H}\Sing(\Gamma_j)\leq m-7.$$

Let $\mathcal R_j$ be any connected component of $\Reg(\Gamma_j)$
and set $\Sigma_j=\overline{\mathcal R_j}$.
We first bound the area of $\Sigma_j$.
Write $\nu$ for the outward normal
and $A$ for the second fundamental form on $\mathcal R_j$. The first
variation of $\mathcal J_j$ gives 
$$\mathbf H=h_j.$$
Since $\mathcal R_j\subset M\setminus K$ for $j$ large enough,
 we have
$\biRic_M\geq m-1$ on $\mathcal R_j$.
Let $\vartheta_j$ denote the least eigenvalue of the Ricci tensor of
$\mathcal R_j$ with its induced metric. At each point, choose a unit
tangent vector $e$ such that 
$$\Ric_{\mathcal R_j}(e,e)=\vartheta_j.$$
The Gauss equation
and the Cauchy--Schwarz inequality on $e^\perp\subset T\mathcal R_j$ give
$$
\begin{aligned}
 &\Ric_M(\nu,\nu)+|A|^2+\vartheta_j\\
 &\quad\geq \biRic_M(\nu,e)+h_j A(e,e)
       +\frac{(h_j-A(e,e))^2}{m-1}\\
 &\quad=\biRic_M(\nu,e)
       +\frac{\left(A(e,e)+\frac{m-3}{2}h_j\right)^2}{m-1}
       -\frac{m-5}{4}h_j^2\\
 &\quad\geq m-1-\frac{m-5}{4}h_j^2.
\end{aligned}
$$
The stability inequality for $\mathcal J_j$, together with
\eqref{eq:annular-prescription-bounds}, now gives
$$
\begin{aligned}
 &\int_{\mathcal R_j} (|\nabla\varphi|^2+\vartheta_j\varphi^2)\,d\mathcal H_g^m\\
 &\quad\geq\int_{\mathcal R_j}
  \bigl(\Ric_M(\nu,\nu)+|A|^2+\vartheta_j+\partial_\nu h_j\bigr)
  \varphi^2\,d\mathcal H_g^m\\
 &\quad\geq (m-1-\delta_j)\int_{\mathcal R_j}\varphi^2\,d\mathcal H_g^m
\end{aligned}
$$
for every $\varphi\in C_c^\infty(\mathcal R_j)$, where
$\delta_j=(9(m-5)+1)/j^2$.

By Lemma~\ref{lem:ampi-components},
$\Sigma_j$ is a compact AM--PI space with connected regular part
$\mathcal R_j$. Its singular set has Hausdorff dimension at most $m-7$,
so we have
$$\mathcal H_g^{m-3}(\Sigma_j\setminus\mathcal R_j)=0.$$
The Ricci tensor of $\mathcal R_j$ is smooth, so its least
eigenvalue $\vartheta_j$ is locally Lipschitz. By the Gauss equation,
$$
 m\vartheta_j\leq\operatorname{Scal}_{\mathcal R_j}
 =\operatorname{Scal}_M-2\Ric_M(\nu,\nu)+h_j^2-|A|^2.
$$
The ambient curvature and $h_j$ are bounded on the compact set
$\Sigma_j$, so $\vartheta_j$ is bounded above.
Theorem~\ref{thm:singularbishop} therefore gives
$$\mathcal H_g^m(\Sigma_j)\leq(1-\delta_j/(m-1))^{-m/2}\omega_m$$
for all sufficiently large $j$ and every such component. Hence we have 
\begin{equation}\label{eq:sharpbubblearea}
 \limsup_{j\to\infty}\mathcal H_g^m(\Sigma_j)\leq\omega_m.
\end{equation}

We now choose a separating component and apply the radial volume
comparison. Fix a unit-speed geodesic ray $\gamma$ from $p$.
Since $\mathcal H_g^{m-1}(\Sing(\Gamma_j))=0$,
Lemma~\ref{lem:separating-component} allows us to choose $\mathcal R_j$
so that its closure $\Sigma_j$ separates $p$ from a point of
$\gamma((8r_j,\infty))$. Both $B_{r_j}(p)$ and $\gamma((8r_j,\infty))$
are connected and disjoint from $\Gamma_j$. They therefore lie in
different components of $M\setminus\Sigma_j$. Since $M$ has one end,
the component containing $B_{r_j}(p)$ is bounded.

Fix $R>0$.
Let $b$ be the Busemann function of $\gamma$, given by
$$
 b(x)=\lim_{t\to\infty}b_t(x),\,\,\text{where }
 b_t(x)=t-d(x,\gamma(t)).
$$
Choose $j$ with $r_j>R$. For $t>8r_j$, every minimizing segment from
$\gamma(t)$ to $B_R(p)$ meets $\Sigma_j$, and
$d(\gamma(t),\Sigma_j)\geq t-8r_j$. On $B_R(p)$ the distance to
$\gamma(t)$ is at most $t+R$ and ranges over an interval of length
$\operatorname{osc}_{B_R(p)}b_t$. Integrating \eqref{eq:abstractradial}
over this interval by the coarea formula gives
$$
 \Vol B_R(p)\leq
 \operatorname{osc}_{B_R(p)} b_t \cdot
 \left(\frac{t+R}{t-8r_j}\right)^m\mathcal H_g^m(\Sigma_j).
$$
The functions $b_t$ converge uniformly to $b$ on $\overline{B_R(p)}$.
Letting first $t\to\infty$ and then $j\to\infty$,  from
\eqref{eq:sharpbubblearea} we obtain
\begin{equation}\label{eq:busemann-volume}
 \Vol B_R(p)\leq\omega_m\operatorname{osc}_{B_R(p)} b.
\end{equation}
Since $b$ is $1$-Lipschitz, its oscillation on $B_R(p)$ is at most
$2R$. This proves the global bound \eqref{eq:globalsharp}.

Now we have that $M$ has one end and linear volume growth.
By \cite[Theorem~19 and Corollary~23]{SormaniBusemann},
$b$ is bounded below. Since $\gamma(0)=p$, we have $b(p)=0$ and hence
$b\leq R$ on $B_R(p)$. Thus \eqref{eq:busemann-volume} gives
\begin{equation}\label{eq:one-ended-volume}
 \Vol B_R(p)\leq \omega_m\bigl(R-\inf_M b\bigr).
\end{equation}
The limit
$$
 \Theta^1(\infty)=\lim_{R\to\infty}\frac{\Vol B_R(p)}R
$$
exists by \cite[Remark~2.2]{ZhouZhu}.
Dividing \eqref{eq:one-ended-volume} by $R$
and letting $R\to\infty$ gives $\Theta^1(\infty)\leq\omega_m$.
The inclusions $B_{R-a}(p)\subset B_R(p')\subset B_{R+a}(p)$ for
$a=d(p,p')$ and $R>a$ give independence of the center.

Finally, suppose that $M$ is nonorientable. Its
orientation double cover $\pi:\widehat M\to M$ satisfies the same
curvature hypotheses with compact set $\pi^{-1}(K)$. Fix a lift
$\widehat p$ of $p$, and let $\tau$ be the deck involution. By path lifting,
$$
 \pi^{-1}(B_R(p))=B_R(\widehat p)\cup B_R(\tau\widehat p).
$$
Write $V(R)=\Vol B_R(p)$, $\widehat V(R)=\Vol B_R(\widehat p)$, and
$a=d(\widehat p,\tau\widehat p)$. The union has volume $2V(R)$ and
lies between $B_R(\widehat p)$ and $B_{R+a}(\widehat p)$. Thus
\begin{equation}\label{eq:covervolume}
 \widehat V(R)\leq 2V(R)\leq\widehat V(R+a).
\end{equation}
The two balls in the union have the same volume, so
$V(R)\leq\widehat V(R)\leq2R\omega_m$ by the orientable case.
Dividing \eqref{eq:covervolume} by $R$ and letting $R\to\infty$ shows
that the volume growth ratio of $M$ is half that of $\widehat M$.
It is therefore independent of $p$ and at most $\omega_m$.

The round cylinder $\mathbb S^m\times\mathbb R$ has volume
growth ratio $2\omega_m$. This proves sharpness of the global
constant and of the two-ended bound. Its quotient by the free isometry
$(x,t)\mapsto(-x,-t)$ has one end. The same double-cover volume
comparison gives volume growth ratio $\omega_m$ on the quotient,
which proves sharpness of the one-ended bound.
\end{proof}

\begin{proof}[Proof of Corollary~\ref{cor:one-ended-isoperimetry}]
Fix $p\in M$. The manifold $M$ has infinite volume, and $I_M$ is
nondecreasing by \cite[Theorem~3.8(1)]{APPS}.
Fix $v>0$ and choose $r>0$ with $\Vol B_r(p)\geq v$.
By the coarea formula, for every $R>r$,
$$
 (R-r)I_M(v)\leq\int_r^R\mathcal P(B_t(p))\,dt
 =\Vol B_R(p)-\Vol B_r(p).
$$
Dividing this inequality by $R$ and letting $R\to\infty$,
Theorem~\ref{thm:main} gives $I_M(v)\leq\Theta^1(\infty)\leq\omega_{n-1}$.
The monotonicity of $I_M$ gives the desired limit bound.
\end{proof}

\appendix
\section{Compact AM--PI spaces}\label{sec:ampi}

We begin with the unweighted form of the AM--PI condition used in
\cite[Section~2.1]{BHHShiZhu} and \cite[Definition~3.1]{BiZhuDescent}.

\begin{definition}\label{def:ampi}
A complete metric measure space 
$$(X,\mathbf d,\boldsymbol{\mu})$$ is an \emph{$m$-dimensional AM--PI
space} if it has a decomposition 
$$X=\mathcal R\sqcup \mathcal S$$ with the
following properties.
\begin{enumerate}[label=\textup{(\roman*)}]
\item The set $\mathcal S$ is closed and $\boldsymbol{\mu}(\mathcal S)=0$.
\item The length metric $\mathbf d_\ell$, induced from $\mathbf d$
within connected components of $\mathcal R$, induces the subspace
topology on $\mathcal R$. Moreover, there is a smooth
Riemannian $m$-manifold $(M,g)$ and a measure-preserving isometry
$$
 (\mathcal R,\mathbf d_\ell,\boldsymbol{\mu}|_{\mathcal R})
 \cong (M,d_g,\mathrm{dvol}_g).
$$
\item The measure $\boldsymbol{\mu}$ is locally Ahlfors $m$-regular, that is,
$$
 C_K^{-1}r^m\leq\boldsymbol{\mu}(B_r(x))\leq C_K r^m.
$$
\item The space $X$ supports the local $(1,2)$-Poincar\'e inequality
\begin{equation}\label{eq:ampi-poincare}
 \fint_B |v-v_B|\,d\boldsymbol{\mu}
 \leq C_K r\left(\fint_{B_{\lambda_K r}(x)\cap\mathcal R}
 |\nabla_g v|^2\,d\boldsymbol{\mu}\right)^{1/2}
\end{equation}
for every $v\in C_c^\infty(\mathcal R)$, where $B=B_r(x)$ and
$v_B=\fint_B v\,d\boldsymbol{\mu}$.
\end{enumerate}
For each compact $K\subset X$, there are constants $C_K,\lambda_K\geq1$
and $r_K>0$ such that \textup{(iii)} and \textup{(iv)} hold for every
$x\in K$ and $0<r<r_K$.
\end{definition}

In the following, let $(X,\mathbf d,\boldsymbol{\mu})$ be a compact
AM--PI space of dimension $m\geq3$ with decomposition
$X=\mathcal R\sqcup\mathcal S$. We assume further that
$$
 \mathcal H^{m-2}(\mathcal S)=0.
$$
We denote by $W^{1,2}(X)$ the completion of $C_c^\infty(\mathcal R)$
in the norm
$$
 \|v\|_{W^{1,2}(X)}^2=\int_X v^2\,d\boldsymbol{\mu}+\int_{\mathcal R}|\nabla v|^2\,d\boldsymbol{\mu}.
$$

\begin{lemma}\label{lem:ampi-constant}
The constant function $1$ belongs to $W^{1,2}(X)$.
\end{lemma}

\begin{proof}
Fix $\varepsilon>0$ sufficiently small. Since $\mathcal S$ is compact and
$\mathcal H^{m-2}(\mathcal S)=0$, it has a finite cover by balls
$B_i=B_{r_i}(x_i)$ with $r_i<\varepsilon$ and $\sum_i r_i^{m-2}<\varepsilon$.
Let $\eta$ be the minimum of $1$ and the functions $\mathbf d(\cdot,B_i)/r_i$.
It vanishes on $\bigcup_i B_i$ and equals $1$ outside $\bigcup_i B_{2r_i}(x_i)$.
Each distance function is $1$-Lipschitz, so the upper Ahlfors bound gives
$$
 \int_X(1-\eta)^2\,d\boldsymbol{\mu}
 +\int_{\mathcal R}|\nabla\eta|^2\,d\boldsymbol{\mu}
 \leq \sum_i(1+r_i^{-2})\boldsymbol{\mu}(B_{2r_i}(x_i))
 \leq C\varepsilon.
$$
Since $\eta$ is Lipschitz with compact support in $\mathcal R$, it can be
approximated in the Sobolev norm by functions in $C_c^\infty(\mathcal R)$.
Letting $\varepsilon\downarrow0$ proves the lemma.
\end{proof}

\begin{lemma}\label{lem:ampi-weight}
Let $V$ be locally Lipschitz and bounded above on $\mathcal R$.
Suppose that for some $\sigma>0$ and every $v\in C_c^\infty(\mathcal R)$,
\begin{equation}\label{eq:ampi-spectral-bound}
 \int_{\mathcal R}(|\nabla v|^2+Vv^2)\,d\boldsymbol{\mu}
 \geq\sigma\int_{\mathcal R}v^2\,d\boldsymbol{\mu}.
\end{equation}
For every $\gamma>1$, there exists $u\in W^{1,2}(X)$ with
$\essinf_Xu=1$ such that
$$
 -\gamma\Delta_{\mathcal R}u+Vu\geq0\,\,\text{on }\mathcal R.
$$
Moreover, $u\in C^{2,\beta}_{\rm loc}(\mathcal R)$ for every $0<\beta<1$.
\end{lemma}

\begin{proof}
For $v,w\in C_c^\infty(\mathcal R)$, set
$$
 Q(v,w)=\int_{\mathcal R}
   \bigl(\gamma\langle\nabla v,\nabla w\rangle+Vvw\bigr)\,d\boldsymbol{\mu}.
$$
Writing $V_\pm=\max\{\pm V,0\}$, we obtain from \eqref{eq:ampi-spectral-bound}
$$
\begin{aligned}
 \int_{\mathcal R}V_-v^2\,d\boldsymbol{\mu}
 &\leq \|\nabla v\|_2^2+\|V_+\|_\infty\|v\|_2^2\\
 \text{and}\,\,Q(v,v)&\geq(\gamma-1)\|\nabla v\|_2^2+\sigma\|v\|_2^2.
\end{aligned}
$$
The first estimate and the Cauchy--Schwarz inequality extend $Q$
continuously to $W^{1,2}(X)$. The second estimate gives coercivity.
Since $V_+\in L^2(X)$, the Lax--Milgram theorem gives $w\in W^{1,2}(X)$
satisfying
$$
 Q(w,v)=\int_{\mathcal R}V_+v\,d\boldsymbol{\mu}
 \,\,\text{for every }v\in W^{1,2}(X).
$$
By Lemma~\ref{lem:ampi-constant} and the Sobolev truncation,
$v=(1-w)_+$ belongs to $W^{1,2}(X)$. Testing with $v$ gives
$$
 Q(v,v)=Q(1-w,v)=-\int_{\mathcal R}V_-v\,d\boldsymbol{\mu}\leq0.
$$
The coercivity yields $v=0$, so $w\geq1$ almost everywhere.
On $\mathcal R$, $w$ is a weak solution of $-\gamma\Delta_{\mathcal R}w+Vw=V_+$.
In local coordinates, the interior estimate for weak solutions gives
$w\in W^{2,2}_{\rm loc}(\mathcal R)$ \cite[Theorem~8.8]{GilbargTrudinger}.
Since the metric is smooth and $V$ and $V_+$ are locally Lipschitz,
\cite[Theorem~9.19]{GilbargTrudinger} with $k=1$ gives
$w\in C^{2,\beta}_{\rm loc}(\mathcal R)$ for every $0<\beta<1$.
As $V_+\geq0$, the function $u=w/(\essinf_X w)$ has the desired properties.
\end{proof}

\begin{lemma}\label{lem:ampi-components}
Let $\mathcal R_0$ be a connected component of $\mathcal R$ and let
$X_0$ be its closure in $X$. Then $X_0$ is open in $X$, and
$(X_0,\mathbf d,\boldsymbol{\mu})$ is a compact AM--PI space with
regular part $\mathcal R_0$ and singular set $X_0\cap\mathcal S$.
\end{lemma}

\begin{proof}
By Lemma~\ref{lem:ampi-constant}, choose
$\{v_i\}_{i=1}^\infty\subset C_c^\infty(\mathcal R)$ converging to $1$
in $W^{1,2}(X)$. Since $\mathcal R_0$ is open and closed in $\mathcal R$,
$v_i\chi_{\mathcal R_0}$ belongs to $C_c^\infty(\mathcal R)$.
It converges to $\chi_{\mathcal R_0}$ in $L^2(X)$, and its gradient
tends to zero in $L^2(\mathcal R)$. Passing to the limit in the
Poincar\'e inequality \eqref{eq:ampi-poincare}, we see that
$\chi_{\mathcal R_0}$ is constant almost everywhere on each sufficiently
small ball. By the lower Ahlfors bound, every nonempty open set has
positive measure. Hence such a ball cannot meet both $\mathcal R_0$ and
$X\setminus X_0$. Every ball centered in $X_0$ meets $\mathcal R_0$,
so $X_0$ is open.

Since $X_0$ is compact and open, $B_{\lambda_X r}(x)\subset X_0$
for every $x\in X_0$ and all sufficiently small $r>0$. The Ahlfors
bounds therefore restrict to $X_0$. Extending each
$v\in C_c^\infty(\mathcal R_0)$ by zero on $\mathcal R\setminus\mathcal R_0$
gives the Poincar\'e inequality on $X_0$. Finally,
$X_0\cap\mathcal R=\mathcal R_0$ inherits its Riemannian length and
volume, and $X_0\cap\mathcal S$ is closed and has measure zero.
\end{proof}

We finish by showing that if $E$ locally minimizes perimeter with
prescribed mean curvature and its perimeter measure has compact
support $\Gamma$, then $\Gamma$ is an AM--PI space.

\begin{proposition}\label{prop:ampi-pmc}
Let $U$ be a relatively compact open subset of a smooth Riemannian
manifold $(M^{m+1},g)$ with $m\geq7$, and let $h\in C^\infty(U)$.
Let $E$ be a Caccioppoli set whose perimeter measure has nonempty
compact support $\Gamma=\spt|D\chi_E|\subset U$. Suppose that
\begin{equation}\label{eq:ampi-local-minimality}
 \mathcal P(E;V)-\int_{E\cap V}h\,dV_M
 \leq \mathcal P(F;V)-\int_{F\cap V}h\,dV_M
\end{equation}
for every open $V\Subset U$ and every Caccioppoli set $F$ with
$E\Delta F\Subset V$. With the restricted ambient distance and
the measure $\mathcal H_g^m\llcorner\Gamma$, the set $\Gamma$ is a
compact AM--PI space. Its regular part $\Reg(\Gamma)$ consists of the
points near which $\Gamma$ is a smooth embedded hypersurface. The
singular set $\Sing(\Gamma)=\Gamma\setminus\Reg(\Gamma)$ has
Hausdorff dimension at most $m-7$.
\end{proposition}

To prove the Poincar\'e inequality, we first estimate the area of a
subset of $\Gamma$ in terms of its relative boundary. Write
$$\boldsymbol{\mu}=\mathcal H_g^m\llcorner\Gamma.$$ For $x\in\Gamma$, let $B_r^M(x)$
be the ambient geodesic ball and put $$B_r(x)=\Gamma\cap B_r^M(x).$$

\begin{lemma}\label{lem:pmc-relative-isoperimetry}
Under the hypotheses of Proposition~\ref{prop:ampi-pmc}, there are
$C,\lambda\geq1$ and $r_0>0$ such that for $x\in\Gamma$ and $0<r<r_0$,
$$
 C^{-1}r^m\leq\boldsymbol{\mu}(B_r(x))\leq Cr^m.
$$
Let $D$ be relatively open in $\Gamma$. Suppose its relative boundary
$\partial_\Gamma D=\overline D^{\,\Gamma}\setminus D$ is a smooth
$(m-1)$-dimensional submanifold of $\Reg(\Gamma)$. If
$\boldsymbol{\mu}(D\cap B_r(x))\leq\boldsymbol{\mu}(B_r(x))/2$, then
\begin{equation}\label{eq:pmc-relative-isoperimetry}
 \boldsymbol{\mu}(D\cap B_r(x))\leq Cr\,
 \mathcal H_g^{m-1}(\partial_\Gamma D\cap B_{\lambda r}(x)).
\end{equation}
\end{lemma}

\begin{proof}
Choose an open set $W$ with $\Gamma\subset W\Subset U$. For every
ball $B_s^M(x)\Subset W$ and every Caccioppoli set $F$ with
$E\Delta F\Subset B_s^M(x)$,
\eqref{eq:ampi-local-minimality} gives
\begin{equation}\label{eq:pmc-ball-comparison}
 \mathcal P(E;B_s^M(x))-\mathcal P(F;B_s^M(x))
 \leq \|h\|_{L^\infty(W)}\Vol(E\Delta F)\leq Cs^{m+1}.
\end{equation}
In Euclidean space, the perimeter measure is $m$-dimensional Hausdorff
measure restricted to the reduced boundary \cite[Theorem~15.9]{Maggi}.
In local coordinates, the Gauss--Green formula expresses Riemannian
perimeter as an integral over the reduced boundary with the area density
induced by $g$ on its approximate tangent planes. The area formula
\cite[3.2.20 and 3.2.46]{Federer} therefore gives
$|D\chi_E|=\mathcal H_g^m\llcorner\partial^*E$ in $U$.
From \cite[(4.75)]{Simmons}, we have $\mathcal H_g^m(\Gamma\setminus\partial^*E)=0$, and
so $|D\chi_E|=\boldsymbol{\mu}$. The perimeter bounds \cite[(4.74)]{Simmons} in a
finite cover of $\Gamma$ by coordinate neighborhoods give
$$C^{-1}s^m\leq\boldsymbol{\mu}(B_s(x))\leq Cs^m$$
uniformly for $x\in\Gamma$ and all sufficiently small $s$.

Let $\beta\in(0,1)$ be the constant in
\cite[Theorem~2]{BombieriGiusti}, and put $L=2/\beta$. Fix
$\eta\in(0,1/2)$. We first claim that there are $c_\eta,s_\eta>0$
such that if $x\in\Gamma$, $0<s<s_\eta$, and
$$
 \min\{\boldsymbol{\mu}(D\cap B_s(x)),\boldsymbol{\mu}(B_s(x)\setminus D)\}
 \geq\eta\,\boldsymbol{\mu}(B_s(x)),
$$
then
\begin{equation}\label{eq:pmc-interface}
 \mathcal H_g^{m-1}(\partial_\Gamma D\cap B_{Ls}(x))
 \geq c_\eta s^{m-1}.
\end{equation}
Otherwise, we can choose $\{(x_i,s_i,D_i)\}_{i=1}^\infty$ satisfying the
density condition, with $s_i\downarrow0$ and
$$
 s_i^{1-m}\mathcal H_g^{m-1}(\partial_\Gamma D_i\cap B_{Ls_i}(x_i))\to0.
$$
Dilate the normal coordinates at $x_i$ by $s_i^{-1}$. Denote the rescaled
boundary current of $E$ by $T_i$ and the image of $D_i$ by $\widetilde D_i$.
On each fixed Euclidean ball $\mathbb B_R=B_R^{\mathbb R^{m+1}}(0)$,
the rescaled metrics converge smoothly to the Euclidean metric, and the
error in \eqref{eq:pmc-ball-comparison} is $O(s_iR^{m+1})$ after rescaling.
According to \cite[Propositions~4.6 and~4.7]{Simmons}, after passing to a subsequence,
$T_i$ converges weakly to a Euclidean perimeter-minimizing boundary $T$
and $\|T_i\|$ converges weakly to $\|T\|$ on every such ball.

In the open ball $\mathbb B_L$, set
$$
 P_i=T_i\llcorner\widetilde D_i
 \,\,\text{and}\,\,
 Q_i=T_i-P_i.
$$
Since $\partial T_i=0$ there, we have $\partial Q_i=-\partial P_i$.
By Stokes' theorem, $\partial P_i$ in $\mathbb B_L$ is the current of
integration with multiplicity one over an oriented portion of the
rescaled image of $\partial_\Gamma D_i$. Hence
$$
 \|\partial P_i\|(\mathbb B_L)
 \leq (1+o(1))s_i^{1-m}
       \mathcal H_g^{m-1}(\partial_\Gamma D_i\cap B_{Ls_i}(x_i))\to0.
$$
By the compactness theorem for integral currents \cite[4.2.17]{Federer},
we may pass to a further subsequence on which $P_i\rightharpoonup P$ and
$Q_i\rightharpoonup Q$ in $\mathbb B_L$. There $P+Q=T$ and
$\partial P=\partial Q=0$. Since $\|P_i\|+\|Q_i\|=\|T_i\|\rightharpoonup
\|T\|$, lower semicontinuity of mass gives $\|P\|+\|Q\|\leq\|T\|$.
The reverse inequality follows from $P+Q=T$. Equality therefore holds,
and the individual mass measures converge weakly to $\|P\|$ and $\|Q\|$.
By the density condition and
the lower area bound at $x_i$, both $\|P_i\|(\mathbb B_1)$ and
$\|Q_i\|(\mathbb B_1)$ are at least $c\eta$ for some $c>0$ independent of $i$.
Weak convergence gives the same lower bound for $\|P\|$ and $\|Q\|$
on the closed unit ball, which is contained in $\mathbb B_2$. But
\cite[Theorem~2]{BombieriGiusti} applied to $T=P+Q$ in $\mathbb B_L$
with $\beta L=2$ gives
$$
 \min\{\|P\|(\mathbb B_2),\|Q\|(\mathbb B_2)\}^{(m-1)/m}
 \leq C\|\partial P\|(\mathbb B_L)=0.
$$
This leads to a contradiction and proves \eqref{eq:pmc-interface}.

Now let $\boldsymbol{\mu}(D\cap B_r(x))\leq\boldsymbol{\mu}(B_r(x))/2$. The area bounds give
a constant $K\geq2$ such that $\boldsymbol{\mu}(B_{2s}(z))\leq K\boldsymbol{\mu}(B_s(z))$
for all $z\in\Gamma$ and small $s$. Take $\eta=1/(2K^2)$ and $r$
small enough that $2r<s_\eta$ and the area bounds hold up to radius
$10Lr$. For every $y\in D\cap B_r(x)$,
$B_r(x)\subset B_{2r}(y)\subset B_{4r}(x)$. Hence
$$
 \boldsymbol{\mu}(B_{2r}(y)\setminus D)\geq\tfrac12\boldsymbol{\mu}(B_r(x))
 \geq\eta\,\boldsymbol{\mu}(B_{2r}(y)).
$$
Let $s_y$ be the smallest radius of the form $2^{1-k}r$, $k\geq0$,
for which $\boldsymbol{\mu}(B_{s_y}(y)\setminus D)\geq\eta\boldsymbol{\mu}(B_{s_y}(y))$.
It exists since all sufficiently small balls about $y$ lie in $D$.
At radius $s_y/2$ the fraction occupied by $D$ exceeds $1-\eta$.
Doubling then gives
$$
 \boldsymbol{\mu}(D\cap B_{s_y}(y))\geq (1-\eta)K^{-1}\boldsymbol{\mu}(B_{s_y}(y))
 \geq\eta\,\boldsymbol{\mu}(B_{s_y}(y)).
$$
Thus \eqref{eq:pmc-interface} applies at $y$ with radius $s_y$.

By the $5r$-covering lemma there is a countable set $Y\subset D\cap B_r(x)$
such that the balls $B_{Ls_y}(y)$, $y\in Y$, are pairwise disjoint
and the balls $B_{5Ls_y}(y)$ cover $D\cap B_r(x)$. Since $s_y\leq2r$,
each selected ball $B_{Ls_y}(y)$ lies in $B_{(1+2L)r}(x)$. The upper
area bound on the covering balls and \eqref{eq:pmc-interface} on the
disjoint balls give
$$
 \begin{aligned}
 \boldsymbol{\mu}(D\cap B_r(x))&\leq C\sum_{y\in Y}s_y^m\\
 &\leq Cr\sum_{y\in Y}\mathcal H_g^{m-1}(\partial_\Gamma D\cap B_{Ls_y}(y))\\
 &\leq Cr\,\mathcal H_g^{m-1}(\partial_\Gamma D\cap B_{(1+2L)r}(x)).
 \end{aligned}
$$
This proves \eqref{eq:pmc-relative-isoperimetry} with $\lambda=1+2L$.
\end{proof}

\begin{proof}[Proof of Proposition~\ref{prop:ampi-pmc}]
By \eqref{eq:pmc-ball-comparison}, $E$ almost minimizes perimeter with
error $Cs^{m+1}$ in balls of radius $s$ near $\Gamma$. Since the metric coefficients
are smooth, \cite[Theorem~1.1]{Simmons} gives local $C^{1,\alpha}$
graphs outside a relatively closed subset of $\Gamma$ of Hausdorff
dimension at most $m-7$.
On each graph, the first variation of the functional in
\eqref{eq:ampi-local-minimality} gives $\mathbf H=h$ in the weak sense
with respect to the outward unit normal of $E$. Since $h$ is smooth,
elliptic regularity implies that the graph is smooth. Thus $\Sing(\Gamma)$
is closed in $\Gamma$, has Hausdorff dimension at most $m-7$, and has
$\boldsymbol{\mu}$-measure zero.

On each component of $\Reg(\Gamma)$, the induced metric $g_\Gamma$
gives the same curve lengths as the ambient metric. Then its Riemannian
distance is therefore the length distance on the component and induces
the subspace topology. By the area formula, its volume measure is the
restriction of $\boldsymbol{\mu}$. Lemma~\ref{lem:pmc-relative-isoperimetry} gives
the local Ahlfors $m$-regularity of $\boldsymbol{\mu}$.

To prove \eqref{eq:ampi-poincare}, fix $B=B_r(x)$ with $x\in\Gamma$
and $0<r<r_0/\lambda$ as in Lemma~\ref{lem:pmc-relative-isoperimetry}.
Extend $v\in C_c^\infty(\Reg(\Gamma))$ by zero to $\Gamma$ and choose
a median $a$ of $v$ on $B$. For $t>0$, each of the sets
$\{v>a+t\}$ and $\{v<a-t\}$ meets $B$ in a set of measure at most
$\boldsymbol{\mu}(B)/2$. Since $v$ vanishes near $\Sing(\Gamma)$, Sard's theorem
implies that for almost every $t>0$ their relative boundaries are the
smooth level sets $\{v=a+t\}$ and $\{v=a-t\}$ in $\Reg(\Gamma)$.
Apply \eqref{eq:pmc-relative-isoperimetry} to the two sets separately,
add the estimates, and integrate over $t>0$. The layer-cake formula
on $B$ and the coarea formula on $B_{\lambda r}(x)\cap\Reg(\Gamma)$ give
$$
 \int_B|v-a|\,d\boldsymbol{\mu}
 \leq Cr\int_{B_{\lambda r}(x)\cap\Reg(\Gamma)}
                |\nabla_{g_\Gamma}v|\,d\boldsymbol{\mu} .
$$
By the triangle inequality, $\fint_B|v-v_B|\,d\boldsymbol{\mu}\leq
2\fint_B|v-a|\,d\boldsymbol{\mu}$. The area bounds give
$\boldsymbol{\mu}(B_{\lambda r}(x))\leq C\lambda^m\boldsymbol{\mu}(B)$, so the Cauchy--Schwarz
inequality on $B_{\lambda r}(x)\cap\Reg(\Gamma)$ yields \eqref{eq:ampi-poincare}.
\end{proof}

\end{document}